\documentclass[11pt]{article}
\usepackage{amsmath,amssymb,theorem}
\usepackage{graphicx, enumerate}
\usepackage{subfigure}
\usepackage{psfrag}
\usepackage{placeins}
\usepackage{wrapfig}
\usepackage{booktabs}
\usepackage{chngcntr}

\usepackage{soul}

\counterwithin{figure}{section}
\numberwithin{figure}{section}
\counterwithin{table}{section}
\usepackage{hyperref}
\hypersetup{
  colorlinks,
  citecolor=black,
  filecolor=black,
  linkcolor=black,
  urlcolor=black
}

\newcommand{\floor}[1]{{\left\lfloor #1 \right\rfloor}}
\newcommand{\ds}{\displaystyle}
\newcommand{\mc}{\mathcal}

\newtheorem{thm}{Theorem}[section]
\newtheorem{conj}[thm]{Conjecture}

\newtheorem{claim}[thm]{Claim}
\newtheorem{lem}[thm]{Lemma}
\newtheorem{prop}[thm]{Proposition}
\theorembodyfont{\rmfamily}
\def\pf{\bigskip\noindent {\emph{Proof}.}~~}

\def\dfn#1{{\sl #1}}
\def\es{\emptyset}
\def\less{\setminus}

\def\qed{ \hfill $\blacksquare$}

  \allowdisplaybreaks[4]

\begin{document}
\title{Improved upper bounds for  Gallai-Ramsey numbers of  odd cycles}
\author{
Christian Bosse, Zi-Xia Song\thanks{Corresponding Author.   E-mail address: Zixia.Song@ucf.edu.},  Jingmei Zhang\\
Department  of Mathematics,  University of Central Florida\\
Orlando, FL 32816, USA\\
}

\maketitle
\begin{abstract}

  A Gallai coloring of a complete graph is an edge-coloring such that no triangle   has  all its edges colored differently.  A Gallai $k$-coloring is a Gallai coloring that uses $k$ colors.   Given an integer $k\ge1$ and a graph $H$, the  Gallai-Ramsey number $GR_k(H)$ is the least positive  integer $n$ such that every Gallai $k$-coloring of the complete graph $K_n$   contains a monochromatic copy of $H$.   Gy\'{a}rf\'{a}s,   S\'{a}rk\"{o}zy,  Seb\H{o} and   Selkow proved in 2010  that $GR_k (H) $ is exponential in $k$ if  $H$ is not bipartite,    linear in $k$ if $H$ is bipartite but  not a star, and constant (does not depend on $k$) when $H$ is a star.  Hence,  $GR_k(H)$ is more well-behaved than the classical Ramsey number $R_k(H)$. However, finding exact values of $GR_k (H)$ is  far from trivial, even when $|V(H)|$ is small.  In this paper, we first improve the existing upper bounds for Gallai-Ramsey numbers of odd cycles by  showing that  $GR_k(C_{2n+1}) \le (n\ln n) \cdot 2^k   -(k+1)n+1$ for all  $k \ge 3$ and $n \ge 8$.  
  We then prove that  $GR_k( C_{13})= 6\cdot 2^k+1$ and $GR_k( C_{15})= 7\cdot 2^k+1$ for     all $k\ge1$.    \end{abstract}

\section{Introduction}
 \baselineskip 16pt

All graphs in this paper are finite and simple; that is, they have no loops or parallel edges. Given a graph $G$ and a set $A\subseteq V(G)$,  we use   $|G|$    to denote  the  number of vertices    of $G$, and  $G[A]$ to denote the  subgraph of $G$ obtained from $G$ by deleting all vertices in $V(G)\less A$.  A graph $H$ is an \dfn{induced subgraph} of $G$ if $H=G[A]$ for some $A\subseteq V(G)$.  We use $K_n$, $C_n$ and $P_n$  to denote the  complete graph, cycle, and path on $n$ vertices, respectively.  
For any positive integer $k$, we write  $[k]$ for the set $\{1,2, \ldots, k\}$. We use the convention   ``$A:=$'' to mean that $A$ is defined to be the right-hand side of the relation.

 Given an integer $k \ge 1$ and graphs $H_1,  \ldots, H_k$, the classical Ramsey number $R(H_1,   \ldots, H_k)$   is  the least    integer $n$ such that every $k$-coloring of  the edges of  $K_n$  contains  a monochromatic copy of  $H_i$ in color $i$ for some $i \in [k]$.  Ramsey numbers are notoriously difficult to compute in general. In this paper, we  study Ramsey numbers of graphs in Gallai colorings, where a \dfn{Gallai coloring} is a coloring of the edges of a complete graph without rainbow triangles (that is, a triangle with all its edges colored differently). Gallai colorings naturally arise in several areas including: information theory~\cite{KG}; the study of partially ordered sets, as in Gallai's original paper~\cite{Gallai} (his result   was restated in \cite{Gy} in the terminology of graphs); and the study of perfect graphs~\cite{CEL}. There are now a variety of papers  which consider Ramsey-type problems in Gallai colorings (see, e.g., \cite{chgr, c5c6,GS, exponential, Hall, DylanSong, C9C11, C6C8, C10}).   These works mainly focus on finding various monochromatic subgraphs in such colorings. More information on this topic  can be found in~\cite{FGP, FMO}.  
 
A \dfn{Gallai $k$-coloring} is a Gallai coloring that uses $k$ colors. 
 Given an integer $k \ge 1$ and graphs $H_1,  \ldots, H_k$, the   \dfn{Gallai-Ramsey number} $GR(H_1,  \ldots, H_k)$ is the least integer $n$ such that every Gallai $k$-coloring of $K_n$   contains a monochromatic copy of $H_i$ in color $i$ for some $i \in [k]$. When $H = H_1 = \dots = H_k$, we simply write $GR_k(H)$ and  $R_k(H)$.    Clearly, $GR_k(H) \leq R_k(H)$ for all $k\ge1$ and $GR(H_1, H_2) = R(H_1, H_2)$.    In 2010, 
Gy\'{a}rf\'{a}s,   S\'{a}rk\"{o}zy,  Seb\H{o} and   Selkow~\cite{exponential} proved   the general behavior of $GR_k(H)$.

\begin{thm} [\cite{exponential}]
Let $G$ be a fixed graph  with no isolated vertices 
 and let $k\ge1$ be an integer. Then
$GR_k(G) $ is exponential in $k$ if  $G$ is not bipartite,    linear in $k$ if $G$ is bipartite but  not a star, and constant (does not depend on $k$) when $G$ is a star.			
\end{thm}

It turns out that for some graphs $G$ (e.g., when $G=C_3$),  $GR_k( G)$ behaves nicely, while the order of magnitude  of $R_k(G)$ seems hopelessly difficult to determine.  It is worth noting that  finding exact values of $GR_k (G)$ is  far from trivial, even when $|G|$ is small.
We will utilize the following important structural result of Gallai~\cite{Gallai} on Gallai colorings of complete graphs.

\begin{thm}[\cite{Gallai}]\label{Gallai}
	For any Gallai-coloring $c$ of a complete graph $G$ with $|G| \ge 2$, $V(G)$ can be partitioned into nonempty sets  $V_1, V_2, \dots, V_p$ with $p>1$ so that    at most two colors are used on the edges in $E(G)\less (E(V_1)\cup \cdots\cup  E(V_p))$ and only one color is used on the edges between any fixed pair $(V_i, V_j)$ under $c$, where $E(V_i)$ denotes the set of edges in $G[V_i]$ for all $i\in [p]$. 
\end{thm}

The partition given in Theorem~\ref{Gallai} is  a \dfn{Gallai-partition} of  the complete graph $G$ under  $c$.  Given a Gallai-partition $V_1,  \dots, V_p$ of the complete graph $G$ under $c$, let $v_i\in V_i$ for all $i\in[p]$ and let $\mathcal{R}:=G[\{v_1,  \dots, v_p\}]$. Then $\mathcal{R}$ is  the \dfn{reduced graph} of $G$ corresponding to the given Gallai-partition under $c$. Clearly,  $\mathcal{R}$ is isomorphic to $K_p$.  
By Theorem~\ref{Gallai},  all edges in $\mathcal{R}$ are colored by at most two colors under $c$.  One can see that any monochromatic copy of $H$ in $\mathcal{R}$  will result in a monochromatic copy of $H$ in $G$ under $c$. It is not  surprising  that  Gallai-Ramsey numbers $GR_k( H)$ are related to  the classical Ramsey numbers $R_2(H)$.  Recently,  Fox,  Grinshpun and  Pach~\cite{FGP} posed the following  conjecture on $GR_k( H)$ when $H$ is a complete graph. 

\begin{conj}[\cite{FGP}]\label{Fox} For all  $k\ge1$ and $t\ge3$,
\[
GR_k( K_t) = \begin{cases}
			(R_2(K_t)-1)^{k/2} + 1 & \text{if } k \text{ is even} \\
			(t-1)  (R_2(K_t)-1)^{(k-1)/2} + 1 & \text{if } k \text{ is odd.}
			\end{cases}
\]
\end{conj}

The first case of Conjecture~\ref{Fox} is related to a question of T. A. Brown  from 1983 (see \cite{chgr}):    What is the largest number $f(k)$ of vertices of a complete graph can have such that it is possible to $k$-color its edges so that every  triangle has edges of   exactly two colors? Such a $k$-coloring of the edges of $ K_{f(k)}$ has   neither a  rainbow triangle  nor a monochromatic triangle. Hence, $ GR_k(K_3)=f(k)+1$ and so the first case of Conjecture~\ref{Fox} follows from the main  result of     Chung and Graham~\cite{chgr} from  1983.  A simpler proof of the first case of Conjecture~\ref{Fox} can be found in~\cite{exponential}.    The    case when $t=4$  was recently settled in~\cite{K4}.  Conjecture~\ref{Fox} remains open for all  $t\ge 5$. The next open case, when $t=5$, involves  $R_2(K_5)$.  Angeltveit and   McKay   \cite{K5} recently  proved that $R_2(K_5)\le 48$. It is widely believed that $R_2(K_5)=43$ (see  \cite{K5}).  It is worth noting that Schiermeyer~\cite{GRK5} recently observed that if  $R_2(K_5)=43$, then Conjecture~\ref{Fox} fails  for $K_5$ when $k=3$. 

In this paper, we continue to  study  Gallai-Ramsey numbers of  odd cycles. Using  the same construction given by   Erd\H{o}s, Faudree, Rousseau and Schelp  in 1976 (see Section 2  in \cite{efrs}) for classical Ramsey numbers of odd cycles,   we see that $GR_k(C_{2n+1})\ge n \cdot 2^k + 1$  for all  $k \ge 1$ and $n \ge 2$.  
 General    upper bounds for   $GR_k (C_{2n+1})$ were first studied in \cite{c5c6},  later   improved  in  \cite{Hall},  and then    in \cite{chen}.

\begin{thm}[\cite{Hall}]\label{old bound}
For all $k \ge 1$ and $n \ge 2$, $GR_k(C_{2n+1}) \le (n\ln n)\cdot (2^{k+3} - 3)$.
\end{thm}

\begin{thm}[\cite{chen}]\label{improved0}
For all  $k \ge 2$ and $n \ge 2$, $GR_k(C_{2n+1}) \le (4n+n\ln n) \cdot 2^k   $.
\end{thm}

In this paper we further improve the existing upper bounds for $GR_k (C_{2n+1})$. We prove the following. 

\begin{thm}\label{improved}
For all  $k \ge 3$ and $n \ge 8$, $GR_k(C_{2n+1}) \le (n\ln n) \cdot 2^k   -(k+1)n+1$.
\end{thm}

The exact values of $GR_k(C_{2n+1})$,  when $n$ is small, have attracted more attention recently.  In 2011,  the  lower bound was shown to be the upper bound for $GR_k(C_5)$.

\begin{thm}[\cite{c5c6}]\label{C5}
  For all $k \ge 1$, $GR_k( C_{5}) = 2 \cdot 2^k + 1$.
\end{thm}

Bruce and Song~\cite{DylanSong} considered the next step and determined the exact values of $GR_k( C_7)$ for all integers $k\ge1$.  More recently, Bosse and Song~\cite{C9C11} determined the exact values for $GR_k( C_9)$ and $GR_k( C_{11})$.

\begin{thm}[\cite{DylanSong}]\label{C7}
 For   all $k \ge 1$, $GR_k( C_{7}) = 3 \cdot 2^k + 1$.
\end{thm}
\begin{thm}[\cite{C9C11}]\label{C9C11}
 For $n \in \{4,5\}$ and all $k \ge 1$, $GR_k( C_{2n+1}) = n \cdot 2^k + 1$.
\end{thm}

  Applying the key ideas in \cite{C9C11}, we determine the exact values of Gallai-Ramsey numbers of $C_{13}$ and $C_{15}$. 

\begin{thm}\label{C13C15}
For   $n\in\{6,7\}$ and all  $k \ge 1$, $GR_k( C_{2n+1}) = n \cdot 2^k + 1$.
\end{thm}

  Theorem~\ref{C9C11} and  Theorem~\ref{C13C15}   provide  partial evidence for the first four open cases of the Triple Odd Cycle Conjecture due to Bondy and Erd\H{o}s~\cite{BE}, which states that $R_3(C_{2n+1})=8n+1$ for all integers  $n\ge2$. \L uczak~\cite{Luczak} showed that  $R_3(C_{2n+1}) =8n + o(n)$, as $n\rightarrow \infty$,  and Kohayakawa,
Simonovits and Skokan~\cite{TOCC} announced a proof  in 2005 that the Triple Odd Cycle Conjecture holds when $n$ is sufficiently large.   
We shall make  use of the following   known results in the proof of Theorem~\ref{improved} and Theorem~\ref{C13C15}.

\begin{lem}[\cite{Hall}]\label{min order t}
For $1 \le t \le n$, any Gallai-colored complete graph having a Gallai partition with at least $4\lceil \frac{n}{t}\rceil + 1$ parts each of order at least $t$ contains a monochromatic $C_{2n+1}$.
\end{lem}
\begin{thm}[\cite{Hall}]\label{GR-even}
For all integers $k\ge1$ and $n\ge2$, $GR_k( C_{2n}) \le (n-1)k + 3n$.
\end{thm}

\begin{thm}[\cite{BE}]\label{2n+1} For every  integer $n \ge 2$, 
 $R_2(C_{2n+1}) =4n+1$.
\end{thm}

\begin{thm}[\cite{cycles}]\label{even cycles}
For every  integer  $n \ge 3$, $R_2(C_{2n}) = 3n-1$.
\end{thm}

\begin{prop}[\cite{CS}]\label{C4C6}
$R_2(C_{4}) =6$ and $R_2(C_{6}) =8$.
\end{prop}

\begin{thm}[\cite{Ros1}, \cite{cycles}]\label{R(Cm, Cl)}
For $4 \le m < \ell$ with $m$ even and $\ell$ odd, $R(C_m, C_\ell) = \max \{ \ell -1 + m/2, \, 2m - 1\}.$
\end{thm}

\begin{lem}[\cite{C9C11}]\label{Hamilton}
For all integers  $\ell\ge3$ and $n\ge1$,   let  $n_1, n_2, \ldots, n_\ell$ be positive integers such that $n_i\le n$ for all $i\in[\ell]$ and $n_1+n_2+\cdots+n_\ell\ge2n+1$. Then the complete multipartite graph $K_{n_1, n_2, \ldots, n_\ell}$ has a cycle of length $2n+1$. \end{lem}

We conclude this section by introducing more notation and mentioning a useful result.   For positive integers $n, k$ and $G:=K_{n\cdot 2^k+1}$, let $c$ be any $k$-edge-coloring of $G$ with color classes $E_1, \dots, E_k$. Then $c$ is  \dfn{bad}  if $G$ contains  neither a rainbow $K_3$ nor a monochromatic $C_{2n+1}$ under  $c$. 
For any  $W\subseteq V(G)$ and  any color $i\in [k]$,    $E:=E_i\cap E(G[W])$  is an \dfn{induced matching} in $G[W]$ if  $E$ is a matching in $G[W]$.
For  two disjoint sets $A, B\subseteq V(G)$,  $A$ is \dfn{mc-complete} to $B$ under the coloring $c$ if all the edges between $A$ and $B$  in $G$ are colored the same color under $c$;  and  we simply say $A$ is     \dfn{$j$-complete} to $B$    if all the edges between $A$ and $B$  in $G$ are colored by some color $j\in[k]$    under $c$;  and   $A$ is \dfn{blue-complete}     to $B$    if all the edges between $A$ and $B$  in $G$ are colored  blue  under $c$. 
For convenience, we use   $A \less B$ to denote  $A -B$; and $A \less b$ to denote  $A -\{b\}$ when $B=\{b\}$.  In addition, we will frequently make use of the following result.

\begin{lem}[\cite{C9C11}]\label{n,n+1 Lemma}
Let  $Y, Z\subseteq V(G)$ be   two disjoint sets with $|Y|  \ge   n$ and $ |Z| \ge n$. If   $Y$ is mc-complete, say blue-complete, to  $Z$  under the coloring $c$,  then no vertex  in $ V(G) \setminus (Y \cup Z)$ is blue-complete to $Y \cup Z$ in $G$.  Moreover, if $|Z| \ge n+1$, then $G[Z]$   has no blue edges. Similarly, if  $|Y| \ge n+1$, then   $G[Y]$ has no blue edges.
\end{lem}

\section{Proof of Theorem \ref{improved}}

Let  $n\ge8$ be given as in the statement. For all $k\ge1$, define the function
\[
f(k,n) :=\begin{cases}
2n+1 & \text{if $k=1$}\\
4n+1 & \text{if $k = 2$}\\
(n \ln n)\cdot 2^k  -(k+1)n+1 & \text{if $k \ge 3$.}
\end{cases}
\]   Clearly, $GR_1(C_{2n+1}) \le f(1, n)$ and by Theorem~\ref{2n+1},   $GR_2(C_{2n+1}) \le f(2, n)$.  It suffices to  show that $GR_k(C_{2n+1}) \le f(k, n)$ for all $k \ge 3$.    Let $G:=K_{f(k, n)}$ and let $c:E(G)\rightarrow [k]$ be any Gallai-coloring of $G$. 
Suppose that  $G$ does not  contain any  monochromatic copy of $C_{2n+1}$ under $c$.  Then $c$ is  bad. Among all complete graphs on $f(k, n)$ vertices with a bad Gallai $k$-coloring,  we choose $G$ with $k$ minimum.  	 
Let $X_1,   \ldots, X_k$  be   disjoint subsets of $V(G)$ such that   for each  $i \in [k]$,   $X_i$  (possibly empty) is  mc-complete in color $i$ to $V(G) \less \bigcup_{i=1}^k X_i$.  Choose $X_1, \ldots, X_k$ so that $\sum_{i=1}^k |X_i| \le (k+1)n$ is as large as possible.  Denote $X:= \bigcup_{i=1}^k X_i$.  Then $|X|\le (k+1)n$. 
 We next  prove several  claims.

\begin{claim}\label{X_i bound}
For all $i \in [k]$, $|X_i| \le n-3$.
\end{claim}
\begin{pf}
Suppose $|X_i| \ge n-2$ for some color $i \in [k]$.  We may assume that color $i$ is blue. We next show that 
$|G \less X| \le f(k-1, n) + 3 $. Suppose $|G \less X| \ge f(k-1, n) + 4 $. Let $A$ be a minimal set of vertices  of $G\less X$ such that $G\less (X\cup A)$ has no blue edges. By minimality of $k$, 
$|G \less (X\cup A)| \le f(k-1, n) -1$. Then $|A|\ge 5$ and so $G\less X$ must contain blue edges. Thus $|X_i| \le n-1$, otherwise for any blue edge $uv$ in $G\less X$, we obtain a blue $C_{2n+1}$ by Lemma~\ref{n,n+1 Lemma}. Let  $t:=  n-|X_i|$. Then $t\in\{1,2\}$ because $n-2 \le |X_i| \le n-1$. It follows that $G\less X$ has a blue $H\in\{(2t+1)K_2,  (2t-1)K_2\cup P_3, K_2\cup  2P_{t+1}, tK_2\cup P_{t+2}, P_4\cup (t-1)P_3,  K_2 \cup P_{2t+1}, P_{2t+2}\}$. But then we obtain a blue $C_{2n+1}$ using $n-t$ vertices in $X_i$, all vertices and edges  of   $H$,   and $n+t+1-|H|$ vertices in $V(G)\less (X\cup V(H))$, a contradiction.  This proves that $|G \less X| \le f(k-1, n) + 3 $.  Thus   
\[
|G| = |X|+|G \less X| \le (k+1)n+f(k-1, n) + 3 = \begin{cases}
4n + (4n + 1) + 3, & \text{if $k=3$}\\
(k+1)n + [(n \ln n)\cdot 2^{k-1 } -kn+1] + 3, & \text{if $k \ge 4$}\\
\end{cases}
\]
so that in any case, $|G| < f(k, n)$ for all $k\ge3$ and $n\ge8$, a contradiction. \qed\medskip
\end{pf}

\begin{claim}\label{empty}
 $X_i = \es$ for   some  $i \in [k]$.
\end{claim}
\begin{pf}
Suppose $X_i \ne \es$ for every $i \in [k]$.  By Claim~\ref{X_i bound}, $|X| \le k(n-3)$.  Then 
\[
|G \less X| \ge f(k, n) - k(n-3) =  (n\ln n) \cdot 2^k   -(k+1)n+ 1  - k(n-3) \ge  (n-1)k+3n,
\]
for all $k \ge 3$ and $n \ge 8$. By Theorem~\ref{GR-even}, $G\less X$    contains  a monochromatic  $C_{2n}$,  and thus  $G$ contains  a monochromatic $C_{2n+1}$,  since $X_i \ne \es$ for all $i \in [k]$, a contradiction.
\qed\\
\end{pf}

By Claims \ref{X_i bound} and \ref{empty},  $|X| \le (k-1)(n-3)$.  Consider now a Gallai partition of $G \less X$ with parts $A_1, \ldots, A_p$,  where $p \ge 2$ and $|A_1| \le |A_2| \le \cdots \le |A_p|$.  By Theorem~\ref{2n+1}, $p \le 4n$.  Additionally, let us define the sets
\[
\begin{split}
B &:=  \{a_i \in \{a_1, \ldots, a_{p-1}\} \mid a_ia_p \text{ is colored blue in } \mc{R} \}\\
R &:= \{a_j \in \{a_1, \ldots, a_{p-1}\} \mid a_ja_p \text{ is colored red in } \mc{R} \}
\end{split}
\]
This motivates us to define the related sets in $G$ as $B_G:= \bigcup_{a_i \in B} A_i$ and $R_G:=\bigcup_{a_j\in R} A_j$.  Moreover, we employ the notation $X_r$ to indicate $X_i$ when $i =$ red, and likewise $X_b$ when $i = $ blue.

\begin{claim}\label{BG U RG lower bound}
$|B_G \cup R_G| \ge 2n + 1$.
\end{claim}
\begin{pf}
Suppose  $|B_G \cup R_G| \le 2n$.  Then every vertex in $B_G\cup R_G$ is either    red- or 
blue-complete to $A_p$.  We may assume that $X_1$ is red-complete to $V(G) \less X$ and $X_2$ is blue-complete to $V(G) \less X$. Let $X'_1:=X_1\cup R_G$, $X'_2:=X_2\cup B_G$, and $X'_i:=X_i$ for all $i\in\{3, \ldots, k\}$.  But then  
\[
\left|\bigcup_{i=1}^k X'_i\right|=|X \cup B_G \cup R_G| \le (k-1)(n-3)  + 2n = (k+1) n - 3(k-1) <   (k+1)n,
\]
contrary to  the  choice of $X_1, \ldots, X_k$.  Thus $|B_G \cup R_G| \ge 2n+1$.\qed
\end{pf}\\

\begin{claim}\label{Ap-2ub}

If $|A_p|\le n$, then $|A_{p-2}|\le \floor{\frac{n}{2}} $. 

\end{claim}

\begin{pf}
 Let $q:= \lfloor \frac{n}{2}\rfloor$. Suppose $|A_p|\le n$ but $|A_{p-2}| \ge q+1$.   Then  $|G| - |A_p\cup A_{p-1}\cup A_{p-2}| - |X|\ge f(k, n) - 3n - (n-3)(k-1) \ge 4n$ for all $k \ge 3$ and $n \ge 8$.  Let $B_1$, $B_2$, $B_3$ be a permutation of $A_{p-2}$, $A_{p-1}$, $A_p$ such that $B_2$ is, say, blue-complete to $B_1 \cup B_3$ in $G$.  Let $b_1,\ldots, b_{q + 1} \in B_1$, $b_{q +2} , \ldots, b_{2q+2}  \in B_2$, and $b_{2q + 3}, \ldots, b_{3q +3}  \in B_3$.    Let $A:= V(G) \less (B_1 \cup B_2 \cup B_3 \cup X)$, and define
\[
\begin{split}
B_1^*&:=\{v\in A \mid  v \text{ is blue-complete to }  B_1   \text{ and red-complete to }  B_3\text{  in } G \}\\
B_2^*&:=\{v\in A\mid  v \text{ is blue-complete to   }  B_1\cup B_3 \text{   in } G\}\\ 
B_3^*&:=\{v\in A\mid  v \text{ is red-complete to   }  B_1 \cup B_3 \text{ in } G \}\\ 
B_4^*&:=\{v\in A \mid  v \text{ is red-complete to }  B_1  \text{ and blue-complete to }  B_3\text{  in }  G \}.\\
\end{split}
\]
Then $A=B_1^*\cup B_2^*\cup B_3^* \cup B_4^*$ and so $|A|=|G| - |A_p\cup A_{p-1}\cup A_{p-2}| - |X|\ge 3n$.  Note that $B_1^*,  B_2^*, B_3^*, B_4^* $ are  pairwise disjoint.    
Suppose first that $B_1$ is red-complete to $B_3$ in $G$. By Lemma \ref{n,n+1 Lemma} applied to $B_3^*$ and $ B_1 \cup B_3$,    $|B_3^*| \le n-1$.  Thus  $|B_1^*| + |B_2^*| + |B_4^*| \ge 3n - (n-1) = 2n+1$. 
By symmetry, we may assume that $|B_1^*| +  |B_2^*| \ge n+1$.  We claim that $G[B_1^* \cup B_2^* \cup B_4^*]$ has no blue edges.  Suppose not. Let $uv$ be a blue edge in $G[B_1^* \cup B_2^* \cup B_4^*]$.   Since $|B_1^*| +  |B_2^*| \ge n+1$, let $x_1, \ldots,  x_{q - 1} \in B_1^* \cup B_2^*$ be distinct vertices that are different from $u$ and $v$.  
If $u, v\in B_1^* \cup B_2^*$,  then we find a blue $C_{2n+1}$ with vertices
\[
\begin{cases}
u, v, b_1, b_{q +2}, b_{2q + 3}, b_{q +3}, \ldots, b_{3q +3}, b_{2q +2}, b_2, x_1, b_3, \ldots, x_{q -2}, b_{q}, & \text{if $n$ is even}\\
u, v, b_1, b_{q +2}, b_{2q + 3}, b_{q +3}, \ldots, b_{3q +3}, b_{2q +2}, b_2, x_1, b_3, \ldots, x_{q -2}, b_{q}, x_{q -1}, b_{q+1}, & \text{if $n$ is odd} 
\end{cases}
\]
a contradiction.   Thus we may assume that $v\in B_4^*$.   
  If $u \in B_1^*\cup B_2^*$,  then we find a blue $C_{2n+1}$ with vertices
  \[
  \begin{cases}
   u, v, b_{2q + 3}, b_{q +2}, b_{2q + 4}, \ldots, b_{3q +3}, b_{2q +2}, b_1, x_1, \ldots, x_{q-2}, b_{q - 1}, & \text{if $n$ is even}\\
    u, v, b_{2q + 3}, b_{q +2}, b_{2q + 4}, \ldots, b_{3q +3}, b_{2q +2}, b_1, x_1, \ldots, b_{q - 1}, x_{q- 1}, b_{q}, & \text{if $n$ is odd}
  \end{cases}
  \]
 a contradiction.  Thus  $u, v \in B_4^*$. But then we obtain  a blue $C_{2n+1}$ with vertices 
\[
\begin{cases}
u, v, b_{2q + 3}, b_{q +2}, b_1, x_1, b_2, \ldots, x_{q - 1}, b_{q}, b_{q + 3}, b_{2q + 4}, b_{q + 4}, \ldots, b_{2q + 1}, b_{3q + 2}, &\text{if $n$ is even}\\
 u, v, b_{2q + 3}, b_{q +2}, b_1, x_1, b_2, \ldots, x_{q - 1}, b_{q}, b_{q + 3}, b_{2q + 4}, b_{q + 4}, \ldots, b_{2q + 2}, b_{3q + 3}, & \text{if $n$ is odd}
\end{cases}
\] 
a contradiction.  This proves that  $G[B_1^* \cup B_2^* \cup B_4^*]$ contains no blue edges. 

Since $|B_1^*| + |B_2^*| + |B_4^*| \ge 2n+1$ and $|A_p|\le n$, by Lemma \ref{Hamilton}, $G[B_1^* \cup B_2^* \cup B_4^*]$ has  a red $C_{2n+1}$, a contradiction.  Thus   $B_1$ must be  blue-complete to $B_3$.  Then $|B_1 \cup B_2 \cup B_3|\le 2n$, else we obtain a blue $C_{2n+1}$ in $G[B_1 \cup B_2 \cup B_3]$. By Lemma \ref{n,n+1 Lemma} applied to $B_2 \cup B_2^*$ and $B_1 \cup B_3$,  we see that  $|B_2^*| \le q - 1$.  
If   $|B_1^*| \ge q$, let $x_1, \ldots, x_{q} \in B_1^*$ be distinct vertices. Then we find a blue $C_{2n+1}$ with vertices  
\[
\begin{cases}
b_1, b_{q + 2}, b_{2q + 3}, b_{q + 3}, \ldots, b_{3q + 3}, b_2, x_1, \ldots, b_{q}, x_{q -1}, & \text{if $n$ is even}\\
b_1, b_{q + 2}, b_{2q + 3}, b_{q + 3}, \ldots, b_{3q + 3}, b_2, x_1, \ldots, b_{q}, x_{q -1}, b_{q + 1}, x_{q}, & \text{if $n$ is odd}
\end{cases}
\]
a contradiction.  Thus $|B_1^*| \le  q - 1$, and similarly, $|B_4^*| \le q - 1$. Therefore, 
\[
\begin{split}
|B_3^*| &= |G| - |X|-|B_1 \cup B_2 \cup B_3| - |B_1^*| - |B_2^*| -  |B_4^*|\\
		&\ge f(k, n)-(k-1)(n-3)-2n - (q - 1 )-  (q - 1 )-(q - 1 )  \\
		&\ge 2n+1.
\end{split}
\]
By Lemma \ref{n,n+1 Lemma} applied to $B_3^*$ and  $B_1 \cup B_3$, $G[B_3^*]$ contains no red edges.  But then by Lemma \ref{Hamilton} and the fact that $|A_p|\le n$ and $|B_3^*| \ge 2n+1$,  $G[B_3^*]$ must contain  a blue $C_{2n+1}$, a contradiction.\qed\medskip
\end{pf}

\begin{claim}
$|A_p| \ge n+1$.
\end{claim}
\begin{pf}
Suppose $|A_p| \le n$. Let $r_i:=|\{j\in[p]: |A_j|\ge i\}|$.  Then $ |G \less X|  = \sum_{i=1}^n r_i $. Let $q:=|A_{p-2}|$. By  Lemma \ref{min order t} and Claim \ref{Ap-2ub},
 \[
\begin{split}
|G| 	 &= |X| +   (|A_p|- q) + (|A_{p-1}| - q)+ \sum_{i=1}^{q}  r_i  \\
& \le (k-1)(n-3) + (2n -2q)+  \sum_{i=1}^{q}  4 \left\lceil \frac{n}{i} \right\rceil  \\
&\le (k-1)(n-3) + (2n-2q)+\sum_{i=1}^{q}  4 \left( \frac{n}{i}  + 1 \right) \\
&= (k-1)(n-3) + 2n+2q+4n\sum_{i=1}^{q}    \frac{1}{i}   \\
 &\le\begin{cases}
 (k-1)(n-3)  + 2n+8+4n\displaystyle\sum_{i=1}^{4}    \frac{1}{i},  & n\in\{8,9\},\,\,  q=\left\lfloor\frac{n}{2}\right\rfloor=4\\[1.2em]
(k-1)(n-3) + 2n+10+4n\displaystyle\sum_{i=1}^{5}   \frac{1}{i},		
 &   n=10,\,\,  q=\left\lfloor\frac{n}{2}\right\rfloor=5\\[1.2em]
(k-1)(n-3)+2n+  2\ds \left\lfloor\frac{n}{2}\right \rfloor +4n \left( 1 + \displaystyle\int_1^{\left\lfloor\frac{n}{2}\right\rfloor} \frac{1}{x}  \, dx\right),  &  n\ge11,\,\,  q=\left\lfloor\frac{n}{2}\right\rfloor \\[1.2em]
(k-1)(n-3)+2n + 2\ds \left(\left\lfloor\frac{n}{2}\right \rfloor-1\right)+4n \left( 1 +\displaystyle  \int_1^{\left\lfloor\frac{n}{2}\right\rfloor-1} \frac{1}{x}  \, dx\right),  &    n\ge8,\,\, q\le \left\lfloor\frac{n}{2}\right\rfloor-1.
\end{cases}\\ 
&\le \begin{cases}
 (k-1)(n-3)  + 2n+8+\ds\frac{25n}{3},  & n\in\{8,9\},\,\,  q=\left\lfloor\frac{n}{2}\right\rfloor=4\\[1.2em]
(k-1)(n-3) + 2n+10 + \ds \frac{137n}{15},		
 &   n=10,\,\,  q=\left\lfloor\frac{n}{2}\right\rfloor=5\\[1.2em] 
(k-1)(n-3)+2n+ 2\displaystyle\left\lfloor\frac{n}{2}\right \rfloor+4n \left( 1 + \ln \frac{n}{2}\right) ,  &  n\ge11,\,\,  q=\left\lfloor\frac{n}{2}\right\rfloor\\[1.2em]
(k-1)(n-3)+2n+ 2\left(\displaystyle\left\lfloor\frac{n}{2}\right \rfloor-1\right)+ 4n \left[ 1 +\ln \ds\left(\frac{n}{2}-1\right)\right],  &    n\ge8,\,\, q\le \left\lfloor\frac{n}{2}\right\rfloor-1.
\end{cases}\\
&<f(k, n), 
\end{split}
\]
for all $k \ge 3$ and $n \ge 8$, a contradiction.  \qed
\end{pf}

\begin{claim}\label{BG, RG lower bound}
$|B_G| \ge n+ 1$ and $|R_G| \ge n+1$. Moreover, both $X_r$ and $X_b$ are empty, giving $|X| \le (k-2)(n-3)$.
\end{claim}
\begin{pf}
We may assume that $|B_G| \ge |R_G|$.  By Claim \ref{BG U RG lower bound}, $|B_G| \ge n+1$.  Suppose for a contradiction that $|R_G| \le n$.    By Lemma \ref{n,n+1 Lemma}, $X_b = \es$ and  neither $G[A_p]$ nor $G[B_G]$ has blue edges.  By minimality of $k$, $|A_p| \le f(k-1, n) - 1$ and $|B_G| \le f(k-1, n) -1$.  Note that $|X| > (k-2)(n-3)$, otherwise
\[
|G| = |A_p| + |B_G| + |R_G| + |X| \le 2[f(k-1, n) - 1] + n + (k-2)(n-3) < f(k, n)
\]
for all $k \ge 3$ and $n \ge 8$, a contradiction.
By Claim~\ref{X_i bound}, $X_i \ne \es$ for all $i \in [k]$ other than blue. Thus neither $G[A_p]$ nor $G[B_G]$ has monochromatic $C_{2n}$. By Theorem \ref{GR-even}, $|A_p| \le (k-1)(n-1) + 3n-1$ and $|B_G| \le (k-1)(n-1) + 3n-1$.  But then
\[
|G| = |A_p| + |B_G| + |R_G| + |X| \le 2[(k-1)(n-1) + 3n -1] + n + (k-1)(n-3) < f(k, n)
\]
for all $k\ge 3$ and $n \ge 8$, a contradiction. Thus, $|B_G| \ge n+ 1$ and $|R_G| \ge n+1$. Therefore, Lemma \ref{n,n+1 Lemma} implies $X_r = \es$, and thus we have $|X| \le (k-2)(n-3)$.  \qed\\
\end{pf}

\begin{claim}\label{Ap-2 upper bound}
$|A_{p-2}| \le n$.
\end{claim}
\begin{pf}  Suppose $|A_{p-2}|\ge n+1$. Then $n+1 \le |A_{p-2}| \le |A_{p-1}| \le |A_p|$ and so  $\mc{R}[\{a_{p-2}, a_{p-1},a_p\}]$ is not a monochromatic triangle in $\mc{R}$ (else $G[A_p\cup A_{p-1}\cup A_{p-2}]$ has a monochromatic $C_{2n+1}$).   Let $B_1$, $B_2$, $B_3$ be a permutation of $A_{p-2}$, $A_{p-1}$, $A_p$ such that $B_2$ is, say blue-complete,  to $B_1 \cup B_3$ in $G$. Then $B_1$ must be  red-complete to $B_3$ in $G$.  By Claim~\ref{BG, RG lower bound},  $|X| \le (k-2)(n-3)$.  Let $A:=V(G)\less (B_1\cup B_2\cup B_3\cup X)$.  By Claim~\ref{n,n+1 Lemma} again, $G[B_2]$ has no blue edges, and neither $G[B_1]$ nor $G[B_3]$ has red or blue edges.  By minimality of $k$,  $|B_1|\le f(k-2, n) -1$ and $|B_3|\le f(k-2, n)-1$. Observe that 
\[ 
\begin{split}
 |A\cup B_2|&=|G|-|B_1|-|B_3| - |X|\\
 		&= f(k, n) - 2[f(k-2, n)-1] - (k-2)(n-3)\\
 		&= \begin{cases}
 		(8n \ln n - 4n + 1)- 2(2n) - (n-3), & k=3\\
 		(16n \ln n - 5n + 1)- 2(4n)- 2(n-3), & k=4\\
 		[(n \ln n)\cdot 2^k - (k+1)n + 1] - 2[(n \ln n)\cdot 2^{k-2}   - (k-1)n] - (k-2)(n-3), & k\ge 5.
 		\end{cases}
 \end{split}
 \] 
 In any case, we see that $|A \cup B_2|  \ge f(k-1, n)$.
 By minimality of $k$, $G[A\cup B_2]$ must have blue edges. By Claim~\ref{n,n+1 Lemma},  no vertex in $A$ is red-complete to $B_1\cup B_3$ in $G$,  and no vertex in $A$ is blue-complete to $B_1\cup B_2$ or $B_2\cup B_3$ in $G$. This implies that 
$A$ must be  red-complete to $B_2$ in $G$. It follows that $G[A]$ must contain a blue edge, say $uv$.    
 Let $b_1, \ldots, b_{n-1}\in B_1$, $ b_{n}, \ldots, b_{2n-2}\in B_2$, and $b_{2n-1}\in B_3$.   If $\{u,v\}$   is blue-complete to $ B_1$, then we obtain a blue $C_{2n+1}$ with vertices $b_1, u, v, b_2, b_n, b_{2n-1}, b_{n+1}, b_3, b_{n+2}, \ldots, b_{n-1}, b_{2n-2}$
   in order, a contradiction.  Thus $\{u,v\}$   is not blue-complete to $ B_1$. Similarly, $\{u,v\}$   is not blue-complete to $ B_3$.
Since no vertex in $A$ is red-complete to $B_1\cup B_3$,  we may assume that $u $ is blue-complete to $ B_1$ and $v $ is blue-complete to $ B_3$. But then we obtain a blue $C_{2n+1}$ with vertices $b_1, u, v, b_{2n-1}, b_n, b_2, b_{n+1},\ldots,  b_{n-1}, b_{2n-2}$ in order.  \qed \\\end{pf}

We may assume that $A_{p-1} \subseteq B_G$. By Claim \ref{Ap-2 upper bound}, $|A_{p-2}|\le n$.   By Lemma \ref{Hamilton},   $|R_G| \le 2n$ and $|B_G \less A_{p-1}| \le 2n$. By Claim \ref{BG, RG lower bound}, 
$|X| \le (k-2)(n-3)$. By minimality of $k$,   $|A_p| \le f(k-2, n) -1$.  
Then
\[
\begin{split}
|G| &= |A_p| + |A_{p-1}| + |B_G \less A_{p-1}| + |R_G| + |X| \\
		&\le 2[f(k-2, n) - 1] + 2n  + 2n + (k-2)(n-3)\\
		&= \begin{cases}
		2(2n) + 2n  + 2n + (n-3), & k = 3\\
		2(4n) + 2n + 2n + 2(n-3), & k = 4\\
		2[(n\ln n) \cdot 2^{k-2}   - (k-1)n] + 2n + 2n + (k-2)(n-3), & k \ge 5. 
		\end{cases}\\
\end{split}
\]
 In any case, we see that $|G|<f(k, n)$ for all  $k\ge3$ and  $n \ge 8$, a contradiction.  This completes the proof of Theorem~\ref{improved}. \qed\\

\section{Proof of Theorem \ref{C13C15}}
Let $n\in\{6,7\}$.  It suffices to  show that $GR_k( C_{2n+1}) \le n \cdot 2^k + 1$ for all $k \ge 1$.  This is trivially true for $k=1$. By Theorem \ref{2n+1} and the fact that $GR_2(C_{2n+1})=R_2(C_{2n+1})$, we may assume that $k\ge3$. Let $G:=K_{n\cdot 2^k + 1}$ and let $c:E(G)\rightarrow [k]$ be any Gallai coloring of $G$.  
We next show that  $G$ contains a  monochromatic copy of $C_{2n+1}$ under the coloring $c$. 
	
 Suppose that  $G$ does not  contain any  monochromatic copy of $C_{2n+1}$ under $c$.  Then $c$ is  bad. Among all complete graphs on $n\cdot 2^k+1$ vertices with a bad $k$-edge-coloring,  we choose $G$ with $k$ minimum.  We next  prove a series of  claims.	

\begin{claim}\label{l-vertex}
Let $W\subseteq V(G)$  and let $\ell\ge 3$ be an integer.  Let $x_1, \ldots, x_\ell \in V(G) \setminus W$ such that   $\{x_1, \ldots, x_\ell\}$ is mc-complete, say blue-complete, to  $W$  under $c$. Let    $q\in\{0, 1, \ldots, k-1\} $ be the number of colors, other than blue,  missing on $ G[W]$ under $c$.   
\begin{enumerate}[(i)]

\item   If $\ell \ge n $, then $|W| \le  n \cdot 2^{k-1-q}$. 
\item   If $\ell = n-1$, then $|W| \le  
n \cdot 2^{k-1-q} +2$.
\item If $\ell = n-2$, then $|W| \le (21-2n) \cdot 2^{k-1-q} + (5n - 31)$
\item If $\ell = n-3$, then $|W| \le 11 \cdot 2^{k-1-q} + (n-7)$
\item If $\ell = n-4$, then $n=7$ and $|W| \le  13 \cdot 2^{k-1-q}$. 
\end{enumerate}

\end{claim}
Put another way, Claim~\ref{l-vertex} asserts that, in particular,
\[
|W| \le \begin{cases}
(2n-1) \cdot 2^{k-1-q}+(n-7), & \text{if $\ell \ge 3$}\\
(2n-3) \cdot 2^{k-1-q}+(n-7), & \text{if $\ell \ge 4$}.
\end{cases}
\]

\pf Each statement (i)-(v) is trivially true if  $|W| < \max\{2n+1 -\ell, n+1\}$. Thus, we may assume that  $|W| \ge \max\{2n+1 -\ell, n+1\}$.  Note that $q\le k-1$. If $q=k-1$, then  all the edges of $G[W]$ are colored only blue. Since  $\{x_1,\ldots, x_\ell\}$ is blue-complete to $W$ and   $|W| \ge  \max\{2n+1 -\ell, n+1\}$, we see that   $G[W\cup \{x_1,  \ldots, x_\ell\}]$ contains a   blue $C_{2n+1}$, a contradiction.  Thus   $q \le k-2$.   

First, assume $\ell\ge n$.  Since $|W|\ge n+1$, by Lemma \ref{n,n+1 Lemma}, $G[W]$ contains no blue edges. By minimality of $k$, $|W|\le n \cdot 2^{k-1-q}$, establishing (i).

For the remainder of the proof, we may assume that $\ell\le n-1$, and that $G[W]$ contains at least one blue edge, otherwise $|W|\le n \cdot 2^{k-1-q}$ by minimality of $k$, giving the result.  Let $W^*$ be a minimal set of vertices in $W$ such that $G[W\less W^*]$ has no blue edges. By  minimality of $k$, $|W\less W^*|  \le n \cdot 2^{k-1-q}$. 

Let $P$ be a longest blue path in $G[W]$  with vertices $v_1, \dots, v_{_{|P|}}$ in order, where $|P|\ge2$.  It can be easily checked that if $|P|\ge 2(n-\ell)+2$, or $|P|= 2(n-\ell)+1$ along with a blue edge in $G[W\less V(P)]$, then  $G[W\cup\{x_1, \ldots, x_\ell\}]$ has a blue $C_{2n+1}$, a contradiction.  Thus   $|P|\le 2(n-\ell)+1$. Assume $|P|= 2(n-\ell)+1$. Then $G[W\less \{v_2, \ldots,v_{_{|P|}}\} ]$ has no blue edges. By minimality of $k$, $|W|=|W\less \{v_2, \ldots,v_{_{|P|}}\}|+(|P|-1)\le n \cdot 2^{k-1-q}+(|P|-1)=n \cdot 2^{k-1-q}+2(n-\ell)$, as desired for each $\ell\in \{n-1, n-2, n-3, n-4\}$.  Thus  $  2\le |P|\le 2(n-\ell)$. 
  
 We now consider the case  $\ell = n-1$. Then $|P|= 2$. If $G[W]$ contains three blue edges, say $u_1w_1, u_2w_2, u_3w_3$, such that $u_1, u_2, u_3, w_1, w_2, w_3$ are all distinct, then we obtain a blue $C_{2n+1}$ with vertices 
\[
\begin{cases}
x_1, u_1, w_1, x_2, u_2, w_2, x_3, u_3, w_3, x_4, u_4, x_5, u_5, & \text{if $n=6$}\\
x_1, u_1, w_1, x_2, u_2, w_2, x_3, u_3, w_3, x_4, u_4, x_5, u_5, x_6, u_6, & \text{if $n=7$}
\end{cases}
\] 
in order, where  $u_4, u_5, u_6\in W\less \{u_1, u_2, u_3, w_1, w_2, w_3\}$, a contradiction. Thus $|W^*|\le 2$ because $|P|=2$. Hence,    $|W|=|W\less W^*|+|W^*|\le n \cdot 2^{k-1-q} + 2$.  This establishes (ii). 
 
\begin{figure}
\centering
\includegraphics[scale=.8]{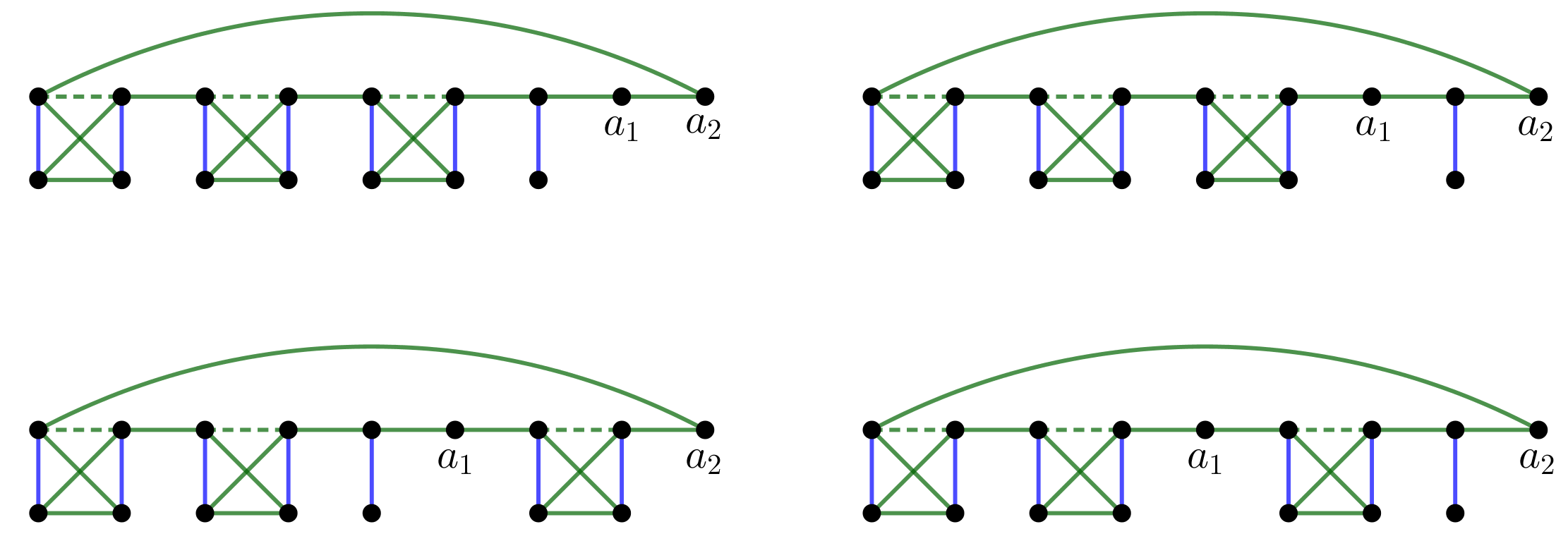}
\caption{The possible ways a green $C_{15}$ arises from a green $C_9$ when $|W| \ge 2m+2$}
\label{C15fromC9}
\end{figure}

 From the above argument,  $\ell  \in\{ n-2, n-3, n-4\}$. Assume first that $|P|=2$. Then   all  the blue edges of $G[W]$ form a matching. Let  $u_1w_1, \ldots, u_mw_m$ be all the blue edges of $G[W]$.  Let 
\[
A:=\begin{cases}
\{u_1, \ldots, u_m\}, &  \text{if $|W|=|\{u_1, \ldots, u_m, w_1, \ldots, w_m\}|=2m$}\\
\{u_1, \ldots, u_m\}\cup \{a\}, & \text{if $|W| - 2m \ge 1$ and $a\in W\less \{u_1, \ldots, u_m, w_1, \ldots, w_m\}$}\\
\{u_1, \ldots, u_m\}\cup \{a_1, a_2\}, & \text{if $n=7$, $|W| - 2m \ge 2$ and $a_1, a_2 \in W\less \{u_1, \ldots, u_m, w_1, \ldots, w_m\}$}\\
\end{cases}
\]
Suppose $|A| \ge (n-3) \cdot 2^{k-1-q}+1$. By Theorem \ref{C7} and Theorem \ref{C9C11}, $G[A]$ has a   monochromatic, say green, $C_{2n-5}$ with $|V(C_{2n-5}) \cap \{u_1, \ldots, u_m\} | \ge 12- n$.  If $n=6$, then we may assume that $E(C_7) = \{u_1u_2, u_2u_3, \ldots, u_6u_7, u_7u_1\}$. Since $G$ has no rainbow triangles under the coloring $c$, then for any $i \in \{1, 3, 5\}$, $\{u_i, w_i\}$ is green-complete to $\{u_{i+1}, w_{i+1}\}$.  Thus we obtain a green $C_{13}$ from the green $C_7$ by replacing the edge $u_iu_{i+1}$ with the path $u_iw_{i+1}w_iu_{i+1}$ for each $i \in \{1,3,5\}$.  If $n=7$, there are four possible ways, up to permutation, $a_1$ and $a_2$ can be arranged on the green $C_9$, as pictured in Figure~\ref{C15fromC9}.  In a manner similar to that for $n=6$, we therefore obtain a green $C_{15}$, a contradiction.
Thus, $|A| \le (n-3) \cdot 2^{k-1-q}  $.  Therefore,
\begin{equation}\label{|P|=2}
\begin{split}
|W|  &= |W \less A| + |A|\\
		&\le \begin{cases}
		2[(n-3) \cdot 2^{k-1-q}], & \text{if $|W| = 2m$}\\
		(7 \cdot 2^{k-1-q} -1)+ 4 \cdot 2^{k-1-q}, & \text{if $|W| \ge 2m+1$}\\
		  (n \cdot 2^{k-1-q} -2)+ (n-3) \cdot 2^{k-1-q}, & \text{if $n=7$ and $|W| \ge 2m+2$}\\
		\end{cases}\\
		&= \begin{cases}
		(n-3) \cdot 2^{k-q}, & \text{if $|W| = 2m$}\\
		(2n-3) \cdot 2^{k-1-q} - 1, & \text{if $|W| \ge 2m+1$}\\
		11 \cdot 2^{k-1-q} -2, & \text{if $n=7$ and $|W| \ge 2m+2$,}
		\end{cases}
\end{split}
\end{equation}
as desired for each $\ell  \in\{ n-2, n-3, n-4\}$. So we may assume that $3\le |P|\le 2(n-\ell)$.

Next suppose $\ell = n-2$.     Then $  |P|\in\{3,4\}$.  Thus  $|W^*| \le 4$, else we obtain a blue $C_{2n+1}$.  Hence,  $|W| = |W \less W^*| + |W^*| \le n \cdot 2^{k-1-q}+ 4$, thus establishing (iii).

\begin{figure}[htbp]
\begin{center}
\includegraphics[scale=0.6]{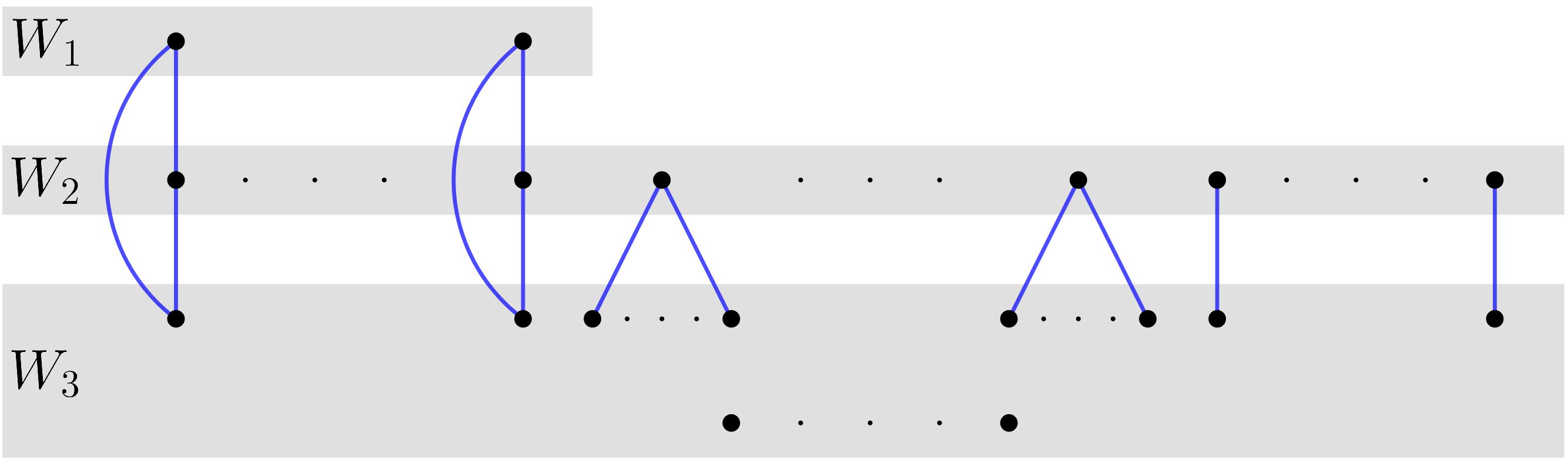}
\end{center}
\caption{Partition of $W$}\label{W1W2W3*}
\end{figure}

By the above arguments, we may now assume that $\ell\in\{n-3, n-4\}$.  Suppose $|P|=3$.  Then each component of the  subgraph of $G[W]$ induced by all its blue edges is isomorphic to a $K_3$, a star, or a $P_2$. Partition $W$ into the sets $W_1$, $W_2$ and $W_3$, described below. 
\[
\begin{split}
W_1: & \text{ Select one vertex from each blue $K_3$}\\
W_2: & \text{ Select one vertex from each blue $K_3$ not in $W_1$, the center vertex}\\
	& \text{ in each blue star, and one vertex from each blue $P_2$}\\
W_3:&=W\less (W_1\cup W_2)
\end{split}
\]
Then  $W = W_1 \cup W_2 \cup W_3$ with $|W |= |W_1| +| W_2|+ | W_3|$.   This partition is illustrated in Figure~\ref{W1W2W3*}.

Note $|W_3| \le n \cdot 2^{k-1-q}$ due to the minimality of $k$.  By an argument similar to the case $|P|=2$, we have $|W_2| \le (n-3) \cdot 2^{k-1-q}$. Therefore, our task is to appropriately bound $W_1$.
Define
\[
A:=\begin{cases}
W_1, & \text{if $n=6$, $\ell = n-3$ and $|W| = 3|W_1|$}\\
W_1\cup \{a\}, & \text{if $n=6$, $\ell=n-3$ and $|W|>3|W_1|$}\\
W_1, & \text{if $n=7$ and $\ell \in \{n-3, n-4\}$},
\end{cases}
\]
where $a\in W$ does not belong to a blue $K_3$.

\begin{figure}[htbp]
\centering
\includegraphics[scale=.8]{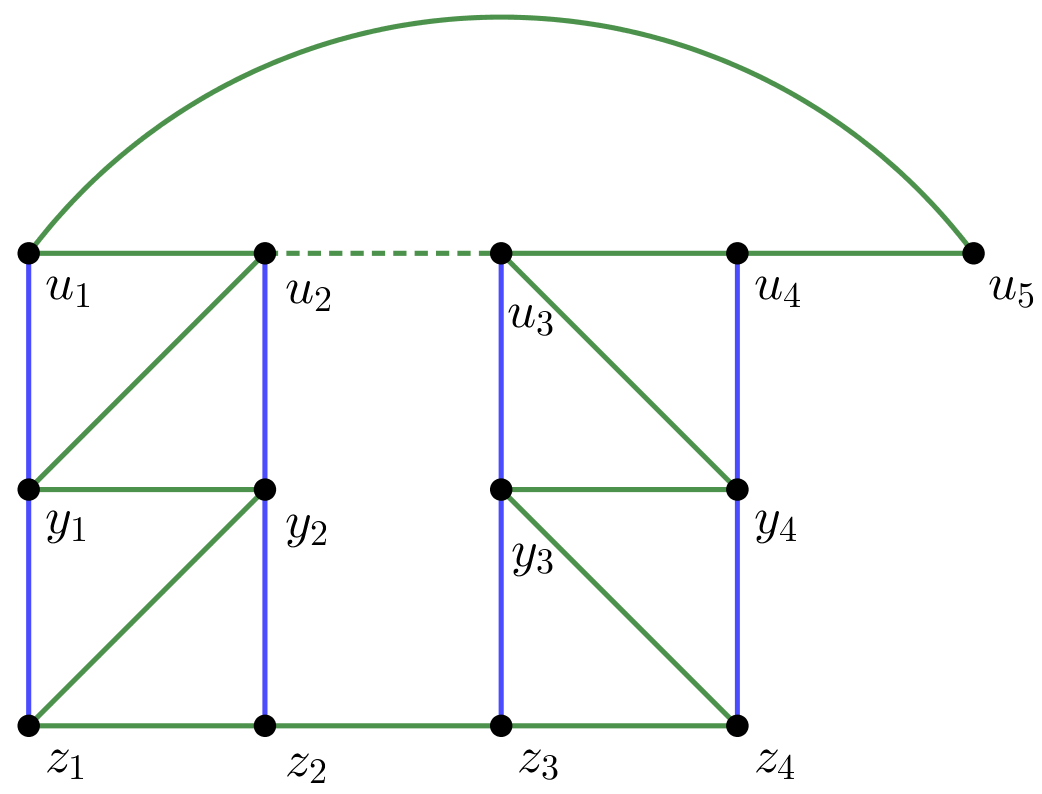}
\caption{A monochromatic $C_{13}$ arising from a monochromatic $C_5$}
\label{C13fromC5}
\end{figure}

We claim that $|A| \le 2 \cdot 2^{k-1-q}$.  First note if $n=7$ and $\ell = n-3$, then $|A| \le 3$ otherwise we find a blue $C_{15}$, giving the result.  Now, suppose $|A| \ge 2 \cdot 2^{k-1-q}+1$.  Then $G[A]$ has a monochromatic, say green, $C_5$ with vertices $u_1, u_2, u_3, u_4, u_5$ (Changed to $u_i$ because $x_i$ already used.)  in order.  We may assume that $a\notin \{u_1, u_2, u_3, u_4\}$.  Enumerate the vertices of the corresponding blue $K_3$'s in $G[W]$ as $u_i, y_i, z_i$ for all $i \in [n-2]$.  Since $G$ has no rainbow triangles under the coloring $c$, then for any $i \in [n-3]$,    $\{u_i, y_i, z_i\}$ is  green-complete to $\{u_{i+1}, y_{i+1}, z_{i+1}\}$.  Additionally, we note $\{u_5, y_5, z_5\}$ is green-complete to $\{u_1, y_1, z_1\}$ when $n=7$. Then we obtain  a green $C_{2n+1}$ with vertices
\[
\begin{cases}
u_1, u_2, y_1, y_2, z_1, z_2, z_3, z_4, y_3, y_4, u_3, u_4, u_5, & \text{if $n=6$}\\
u_1, y_2, y_3, y_4, y_5, y_1, z_2, z_3, z_4, z_5, z_1, z_2, u_3, u_4, u_5, & \text{if $n=7$}
\end{cases}
\]
in order, a contradiction (see Figure~\ref{C13fromC5}).   Thus  $|A| \le 2 \cdot 2^{k-1-q}$.

Therefore,
\begin{align*}
|W| &= |W_1| + |W_2| + |W_3| \\
&\le  \begin{cases} 
3(2 \cdot 2^{k-1-q})&  \text{if $n=6$, $\ell = n-3$ and $|W| = 3|W_1|$}\\
(2 \cdot 2^{k-1-q} - 1) + 3 \cdot 2^{k-1-q} + 6 \cdot 2^{k-1-q} &  \text{if $n=6$, $\ell = n-3$ and $|W|> 3 |W_1|$}\\
 3+3 \cdot 2^{k-1-q} + 6 \cdot 2^{k-1-q} &  \text{if $n=7$ and $\ell = n-3$}\\
2 \cdot 2^{k-1-q}   + 4 \cdot 2^{k-1-q} + 7 \cdot 2^{k-1-q} &  \text{if $n=7$ and $\ell = n-4$}
\end{cases}\\
&=\begin{cases} 
6\cdot 2^{k-1-q} &  \text{if $n=6$, $\ell = n-3$ and $|W| = 3|W_1|$}\\
 11 \cdot 2^{k-1-q} -1 &  \text{if $n=6$, $\ell = n-3$ and $|W|> 3 |W_1|$}\\
9 \cdot 2^{k-1-q} +3 &  \text{if $n=7$ and $\ell = n-3$}\\
 13 \cdot 2^{k-1-q}   &  \text{if $n=7$ and $\ell = n-4$}
\end{cases}
\end{align*}
 as desired. So we may assume that $4\le |P|\le 2(n-\ell)$.  

Now suppose $\ell = n-3$.  Then    $4\le |P| \le 6$.    Assume first that $|W^*| \le 7$. Then 
\[|W| = |W \less W^*| + |W^*| \le n \cdot 2^{k-1-q} + 7 < 11 \cdot 2^{k-1-q} + n-7
\]
because $q \le k-2$ and $k \ge 3$. So we may assume that $|W^*| \ge 8$.   Let $P'$ be a longest blue path in $G[W \less V(P)]$.

Let us first handle the case when $n=6$. Because $|W^*| \ge 8$, we have $|P| \in\{4,5\}$, and
\[
|P'| \le \begin{cases}
3, & \text{if $|P| = 4$}\\
2, & \text{if $|P| = 5$}.
\end{cases}
\]
Moreover, when $|P| = 4$, there is at most one $P'$ such that $|P'| = 3$, otherwise we obtain a blue $C_{13}$.   When $|P|= 4$, it suffices to consider the worst-case scenario, namely when $|P'| = 3$, with vertices $y_1,y_2,y_3$ in order. Define
\[
A:= \begin{cases}
\{v_2, v_3, v_4, y_3\}, & \text{if $|P| = 4$}\\
\{v_2, v_3, v_4, v_5\}, & \text{if $|P| = 5$}
\end{cases}
\]
Then the blue edges of $G[W \less A]$ induce a matching.  Similar to the above case when $|P|=2$, we obtain $|W\less A|\le 9\cdot   2^{k-1-q} - 1$. Hence, 
\[
|W| = |W \less A| + | A| \le 9\cdot   2^{k-1-q} - 1 +3  = 9 \cdot 2^{k-1-q} + 2, 
\]
which is less than the desired bound. 

Now we consider when $n=7$. Again because $|W^*| \ge 8$, we have $|P| = 4$ and $|P'| =2$,  else we obtain a blue $C_{15}$.  Thus  the blue edges in $G[W \less V(P)]$ form an induced matching.  Let  $u_1w_1, \ldots, u_mw_m$ comprise the blue edges of   $G[W \less V(P)]$.   Let 
\[
\begin{split}
A&:= W \less \{v_1, v_2, v_3, v_4, u_1, \ldots, u_m, w_1, \ldots, w_m\}\\
B&:= \begin{cases}
\{v_1, u_1, \ldots, u_m\}\cup A, &  \text{if $|A|\le1$}\\
\{v_1, u_1, \ldots, u_m\}\cup \{a_1, a_2\}, & \text{if $|A|\ge2$, where $a_1, a_2\in A$}
\end{cases}
\end{split}
\]
 By similar reasoning to the case when $|P|=2$, we have  $|B| \le 4 \cdot 2^{k-1-q} $.  Note that when $|A|\ge2$, $|W\less \{v_1,v_2, v_3, u_1, \ldots, u_m\}|\le 7\cdot 2^{k-1-q}$ by minimality of $k$. Therefore,  
\begin{align*}
|W| &= \begin{cases}
 2|B \less (A \cup \{v_1\})| + |A| + |P|  & |A| \le 1  \\
|W\less \{v_1,v_2, v_3, u_1, \ldots, u_m\}|+|B \less \{a_1, a_2\}| + |\{v_2, v_3\}|  & |A| \ge 2,\\
\end{cases}\\
&\le \begin{cases}
 2(4 \cdot 2^{k-1-q}-1)+1+4    & \text{if $|A| \le 1$}\\
7\cdot 2^{k-1-q} + (4\cdot 2^{k-1-q} -2)+2  & \text{if $|A| \ge 2$},\\
\end{cases}
\end{align*} 
  yielding the desired bound because $q \le k-2$ and $k\ge 3$. This establishes case (iv).
 
Finally, we prove case (v), when $\ell = n-4$ and $n=7$.   Assume first that $|W^*| \le 12$. Then 
\[
|W| = |W \less W^*| + |W^*| \le 7 \cdot 2^{k-1-q}+12 \le 13 \cdot 2^{k-1-q}
\]
because $q \le k-2$ and $k \ge 3$. Thus we may assume that $|W^*| \ge 13$. Hence,  $4\le |P| \le 7$, else we obtain a blue $C_{15}$.   Again we will let $P'$ be a longest blue path in $G[W \less V(P)]$.

Let us first handle the cases when $|P| \in \{6, 7\}$. Then there exists a subset $A\subseteq W$ such that $|A|\le5$ and all the blue edges in $G[W\less A]$ form a matching.  By similar reasoning to the case $|P|=2$, we have $|W\less A|\le 11\cdot 2^{k-1-q}-1$, which yields
\[
|W|=|W\less A|+|A|\le 11\cdot 2^{k-1-q}-1+5 \le 13\cdot 2^{k-1-q},
\]   
because $q \le k-2$ and $k\ge 3$.

Now suppose $|P| = 5$.  Except for one case, we may apply identical reasoning as when $|P| \in \{6, 7\}$.  The only case we need to consider is when $|P'|=3$ for possibly many disjoint longest blue paths in $G[W \less V(P)]$.  Apply the partition on $G[W \less \{v_1, v_2\}]$ used to derive the case when $|P|=3$, to obtain corresponding parts $W_1'$, $W_2'$ and $W_3'$ (see Figure~\ref{W1W2W3*}).  By similar reasoning, we find $|W_1'| \le 2 \cdot 2^{k-1-q}$, and $|W_3'| \le 7 \cdot2^{k-1-q}$.  From an argument similar to the case $|P|=2$ used to obtain (\ref{|P|=2}), $|W_2'| \le 4 \cdot 2^{k-1-q}-2$.  Adding the parts together, 
\[
\begin{split}
|W|&=|W_1'|+|W_2'| + |W_3'|+|\{v_1, v_2\}| \le 2\cdot 2^{k-1-q} + (4\cdot 2^{k-1-q} - 2) + 7 \cdot 2^{k-1-q} + 2 = 13 \cdot 2^{k-1-q}
\end{split}
\]
since $q \le k-2$ and $k \ge 3$, as desired.

Thus $|P|=4$.  Then $G[W \less V(P)]$ has at most one blue $P_4$, else we obtain a blue $C_{15}$. It suffices to consider the worst-case scenario when $G[W \less V(P)]$ has exactly one blue $P_4$, with vertices $y_1, y_2, y_3, y_4$ in order.  
Then each component of the  subgraph of  $G[W\less\{v_4, y_4\}]$ induced by all its blue edges is isomorphic to a $K_3$, a star, or a $P_2$.  Define the following sets:
\[
\begin{split}
A_0: & \text{ All vertices $v\in W$ such that $v$ is not incident with any blue edge in $G[W]$}\\
A_1: & \text{ Select one vertex from each blue $K_3$}\\
A_2: & \text{ Select one vertex from each blue $K_3$ not in $W_1$, the center vertex}\\
	& \text{ in each blue star, and one vertex from each blue $P_2$}\\
\end{split}
\]
We next choose $W_1$ and  $W_2$ judiciously. If $G[W\less\{v_4, y_4\}]$ has no blue star, let 
\[
\begin{split}
W_1 &:=A_1\\
W_2 &:=\begin{cases}
A_2\cup A_0,& \text{if $|A_0|\le 1$}\\
A_2\cup \{a_1,a_2\}, & \text{if $|A_0|\ge2$, where $a_1, a_2\in A_0$}
\end{cases} \\
W_3 &:=W\less (A_1\cup A_2)\cup \{v_4, y_4\}.
\end{split}
\]
If on the other hand $G[W\less\{v_4, y_4\}]$ has at least one blue star $S$, with center vertex $x$, and two leaves  $y_1, y_2$, let  
\[
\begin{split}
W_1 &:=A_1\cup \{y_1\}\\
W_2 &:=(A_2\less\{x\})\cup\{y_1, y_2\}\\
W_3 &:=W\less (A_1\cup A_2\cup\{v_4, y_4\}).
\end{split}
\]
  By a similar argument to that given for the case $|P|=3$ (with $\ell=n-4$ and $n=7$), $|W_1|\le 2 \cdot 2^{k-1-q}$, $|W_2|\le 4 \cdot 2^{k-1-q}$ and $|W_3|\le 7 \cdot 2^{k-1-q}$.  Therefore,  
\begin{align*}
|W| &= \begin{cases}
 |W_1| + 2|W_2\less A_0 | +|A_0| +| \{v_4, y_4\}| & |A_0| \le 1 \text{ and  } G[W]  \text{ has no blue star}  \\
 |W_1| + |W_2\less\{a_1, a_2\}|+|\{ v_4, y_4\}|+|W_3| & |A_0| \ge 2 \text{ and  } G[W]  \text{ has no blue star}  \\
 |W_1\less \{y_1\}| +|W_2\less\{y_1, y_2\}|+|\{x, v_4, y_4\}|+|W_3|  &  G[W]  \text{ has a blue star},\\
\end{cases}\\
&\le \begin{cases}
 2 \cdot 2^{k-1-q}+8 \cdot 2^{k-1-q}+1+2    & |A_0| \le 1 \text{ and  } G[W]  \text{ has no blue star}  \\
  2 \cdot 2^{k-1-q}+(4 \cdot 2^{k-1-q} -2)+2  +7\cdot 2^{k-1-q} & |A_0| \ge 2 \text{ and } G[W]  \text{ has no blue star} \\
( 2 \cdot 2^{k-1-q}-1)+(4 \cdot 2^{k-1-q} -2)+3+7\cdot 2^{k-1-q}    &  G[W]  \text{ has a blue star},\\
\end{cases}
\end{align*} 
  yielding the desired bound because $q \le k-2$ and $k\ge 3$. This completes the proof of  Claim \ref{l-vertex}.\qed\\

Let $X_1,   \ldots, X_m$  be a maximum sequence of disjoint subsets of $V(G)$ such that,  for all $j \in [m]$, one of the following holds: 
\begin{enumerate}[(a)]
\item\label{property a}    $1 \le |X_j| \le 3$, and $X_j$ is mc-complete to $V(G) \setminus \bigcup_{i \in [j]} X_i$ under $c$, or 

\item\label{property b}     $4 \le |X_j| \le 6$, and $X_j$ can be partitioned into two non-empty sets $X_{j_1}$ and $X_{j_2}$, where $j_1, j_2\in [k]$ are two distinct colors, such that for each  $t\in\{1,2\}$,   $1\le |X_{j_t}|  \le 3$, $X_{j_t}$ is $j_t$-complete to $V(G) \setminus \bigcup_{i \in [j]} X_i$ but not $j_t$-complete to $X_{j_{3-t}}$, and all the edges between  $X_{j_1}$ and $X_{j_2}$ in $G$ are colored using only the colors $j_1$ and $j_2$.  
\end{enumerate}

  \noindent Note that such a sequence $X_1,   \ldots, X_m$  may not exist. Let $X:= \bigcup_{j \in [m]} X_j$.   For each $x\in X$,  let $c(x)$ be the unique color on the edges between $x$ and $V(G) \less X$ under $c$.  For all $i \in [k]$, let $X_i^*:= \{x \in X : c(x)   \text{ is color } i\}$. Then  $X = \bigcup_{i \in [k]} X_i^*$. It is worth noting that for all $i \in [k]$, $X_i^*$ is possibly empty.   By abusing the notation, we use $X_b^*$, $X_r^*$ and $X_g^*$ to denote $X_i^*$  when $i$  is blue, red or green, respectively.
  
\begin{claim}\label{X*}
For all $i \in [k]$, $|X_i^*| \le 3$. Hence,  $|X|\le 3k$. 
\end{claim}

\pf
Suppose the statement is false.  Then $m \ge 2$.  When choosing $X_1, X_2, \ldots, X_m$, let $j\in[m-1]$ be the largest index such that $|X_p^*\cap (X_1\cup X_2\cup\dots \cup X_{j})|\le 3$ for all colors $p\in[k]$. Then $4\le |X_i^*\cap (X_1\cup X_2\cup\dots \cup X_j\cup X_{j+1})|\le 6$ for some color $i\in[k]$ by the choice of $j$.  Such a color $i$ and an index $j$ exist due to the assumption that the statement of Claim \ref{X*} is false.  Let $A :=X_1\cup X_2\cup  \cdots\cup X_j\cup X_{j+1}$. 
By the choice of   $X_1, X_2, \ldots, X_m$,  there are   at most two colors $i\in[k]$ such that $4\le |X_i^*\cap A|\le 6$. We may assume that such a color $i$ is red or blue. 
Let  $A_b:=\{x\in A: \, c(x) \text{ is color blue} \}$ and $A_r:=\{x\in A: \, c(x) \text{ is color red} \}$.  It suffices to consider the worst-case scenario when $4 \le |A_b| \le 6$ and $4 \le |A_r| \le 6$. Then for any color $p\in[k]$ other than  red and blue,  $|X_p^*\cap A|\le 3$.  Thus by the choice of $j$, $|A\less (A_b \cup A_r)|\le 3(k-2)$. We may assume that $|A_b|\ge|A_r|$. Note that $4\le |A_b|\le6\le n$.  
By Claim \ref{l-vertex}  applied to   $A_b$ and $V(G)\less A$,  we see that 
\begin{align*}
|V(G) \less A| &\le \begin{cases}
(2n-3) \cdot 2^{k-1} +(n-7), & \text{if $|A_b| = 4$}\\
n \cdot 2^{k-1} + (2n-10), & \text{if $|A_b| = 5$}\\
n \cdot 2^{k-1} + (2n-12), & \text{if $|A_b| = 6$}
\end{cases}
\end{align*} 
 But then, 
\begin{align*}
|G| &= |A\less (A_b \cup A_r)| + |A_b| + |A_r|+ |V(G) \setminus A| \\
&\le 3(k-2) + \begin{cases}
 4+4+ [(2n-3) \cdot 2^{k-1} +(n-7)], & \text{if $|A_b| = 4$}\\
5 + 5 + [n \cdot 2^{k-1} + (2n-10)],  & \text{if $|A_b| = 5$}\\
6 +6 + [n \cdot 2^{k-1} + (2n-12)], & \text{if $|A_b| = 6$}
\end{cases}\\
&< n \cdot 2^k + 1
\end{align*}
for all $k \ge 3$ and  $n\in\{6,7\}$, a contradiction.
\qed\\

By Claim \ref{X*}, $|X| \le 3k$.  Let $X'\subseteq X$  be such that for  all $i \in [k]$,  $|X'\cap X^*_i|=1$ when $X^*_i \ne \es$.  Similarly, define $X'' \subseteq X$ such that for all $i \in [k]$, $|X'' \cap (X_i^* \less X')| = 1$ when $X_i^* \less X' \ne \es$.  Finally, let $X''':=X \less (X' \cup X'')$.  
Now consider a  Gallai partition $A_1, \ldots, A_p$ of $G \setminus X$ with $p\ge2$. 
We may  assume that  $1\le |A_1| \le \cdots \le  |A_s|<3 \le  |A_{s+1}| \leq \cdots \leq |A_{p}|$, where $0\le  s\le p$. Let  $\mathcal{R}$ be  the reduced graph of $G\less X$ with vertices $a_1, a_2, \dots, a_{p}$, where $a_i\in A_i$ for all $i\in[p]$. By Theorem~\ref{Gallai}, we may assume that  every  edge  of $\mathcal{R}$ is  colored either red or blue. Note that any monochromatic $C_{2n+1}$ in $\mathcal{R}$ would yield a monochromatic $C_{2n+1}$ in $G$. Thus $\mathcal{R}$ has neither a red nor a blue $C_{2n+1}$.  By Theorem~\ref{2n+1}, $p \le 4n$. Then $|A_p|\ge 2$ because $|G\less X| \ge n\cdot 2^k+1 -3k\ge 8n-8$ and $n\in\{6,7\}$.   

\begin{claim}\label{Ap-8}   

 $|A_{p-8}| \le 2$ and    $|A_{p-4n+12}| \le1$. Moreover, if $|A_{p-7}| \ge 3$, then $|A_{p-4n+16}| \le n-6$.   Similarly, if $|A_{p-4n+13}| \ge 2$, then $p\le 4n-12$.  
\end{claim} 

\pf  
Suppose $|A_{p-8}| \ge 3$ or  $|A_{p-7}| \ge 3$ and $|A_{p-4n+16}| \ge n-5$. By Theorem~\ref{2n+1},  $R_2(C_{2n-7})=4n-15$. We see that   either  $\mc{R}[\{a_{p-8}, a_{p-7}, \ldots, a_p\}]$ has a monochromatic $C_5$   that gives a monochromatic $C_{2n+1}$ in $G$,  or  $\mc{R}[\{a_{p-4n+16}, a_{p-4n+15}, \ldots, a_p\}]$ has a monochromatic $C_{2n-7}$ which  again yields a monochromatic $C_{2n+1}$ in $G$, a contradiction.  
Similarly,  suppose $|A_{p-4n+12}| \ge 2$ or $|A_{p-4n+13}| \ge 2$ and $p\ge 4n-11$ (and so $|A_{p-4n+12}|\ge1$). By Theorem~\ref{2n+1},  $R_2(C_{2n-5})=4n-11$. Thus $\mc{R}[\{a_{p-4n+12}, a_{p-4n+13}, \ldots, a_p\}]$ has a monochromatic $C_{2n-5}$, again yielding a monochromatic $C_{2n+1}$ in $G$, a contradiction. 
\qed \\

\begin{claim}\label{Ap4}
$|A_p| \ge 4$.
\end{claim}

\pf
Suppose $  |A_p| \le 3$.  Then  $ n\cdot 2^k+1 -3k\le |G\less X| \le p|A_p|=12n$ because  $p\le 4n$ and $|X|\le 3k$.  It follows that  $k=3$ and so  $|X|\le 3k=9$ and $|G| = 8n + 1$.  Thus $p\ge 2n+1$ because $|A_p|\le3$.  Let green be the third color. Since $|A_p| \le 3$, we see that  $G$ has no green $C_{2n}$ under the coloring $c$.     We claim that  either $|X^*_r|=0$ or $|X^*_b|=0$.  Suppose $|X^*_r|\ge1$ and  $|X^*_b|\ge1$. Since $G$ has no green $C_{2n}$ and 
\[ |G| - |A_p\cup X| \ge (8n+1) - 3 - 9 \ge 8n-11 > 6n-3 \ge GR_3(C_{2n}),\]
  by Theorem \ref{GR-even},  there is either a red or a blue $C_{2n}$ in $G\less (A_p\cup X)$. Thus $G\less (A_p\cup X^*_g)$ has either a red or a blue $C_{2n+1}$ under $c$, a contradiction. This proves that  either $|X^*_r|=0$ or $|X^*_b|=0$.  We may assume that $|X^*_b|=0$. Then $|X'| \le 2$ and so     $|X|=|X_r^*\cup X_g^*|\le 6$.  By Claim \ref{Ap-8}, $|A_{p-8}| \le 2$ and    $|A_{p-4n+12}| \le1$. 
  If $p\le 4n-6$ or $|X|\le4$, then  
 \[
\begin{split}
|G| 	&= \sum_{i=1}^{p} |A_i| + |X| \le  3 \cdot 8 + 2(4n-20) + (p-4n+12) + |X|\le 8n  < 8n+1,
\end{split}
\]
a contradiction.   Thus $p\ge 4n-5$ and $|X|\ge5$. Since $|X'|\le2$ and $|X|\ge5$, by Claim~\ref{X*},  $|X'|=2$ and  $|X^*_r|\ge2$.  By Theorem \ref{R(Cm, Cl)}, $R(C_{2n-2}, C_{2n+1}) = 4n-5$. It follows that  $\mc{R}[\{a_1, \ldots, a_{4n-5}\}]$ has either  a red $C_{2n-2}$ or a blue $C_{2n+1}$. Since $c$ is bad, we see that $\mc{R}[\{a_1, \ldots, a_{4n-5}\}]$ has a red $C_{2n-2}$. But then $G[V(C_{2n-2})\cup X^*_r\cup \{v\}]$, where $v\in A_p$, has a red $C_{2n+1}$, a contradiction.\qed\\

\begin{claim}\label{Ap-2}

If $|A_p|\le n$, then $|A_{p-2}|\le3$. 

\end{claim}

\pf  Suppose $|A_p|\le n$ but $|A_{p-2}| \ge 4$.  Since $|A_p| \le n$, we have $|G| - |A_p\cup A_{p-1}\cup A_{p-2}| - |X|\ge n\cdot 2^k+1 - 3n  - 3k\ge 5n  -8$.  Let $B_1$, $B_2$, $B_3$ be a permutation of $A_{p-2}$, $A_{p-1}$, $A_p$ such that $B_2$ is, say, blue-complete to $B_1 \cup B_3$ in $G$.  This is possible due to Theorem~\ref{Gallai}.   Let $b_1,\ldots, b_4 \in B_1$, $b_5, \ldots, b_8 \in B_2$, and $b_9, \ldots, b_{12} \in B_3$.    Let $A:= V(G) \less (B_1 \cup B_2 \cup B_3 \cup X)$, and define
\[
\begin{split}
B_1^*&:=\{v\in A \mid  v \text{ is blue-complete to }  B_1   \text{ and red-complete to }  B_3\text{  in } G \}\\
B_2^*&:=\{v\in A\mid  v \text{ is blue-complete to   }  B_1\cup B_3 \text{   in } G\}\\ 
B_3^*&:=\{v\in A\mid  v \text{ is red-complete to   }  B_1 \cup B_3 \text{ in } G \}\\ 
B_4^*&:=\{v\in A \mid  v \text{ is red-complete to }  B_1  \text{ and blue-complete to }  B_3\text{  in }  G \}.\\
\end{split}
\]
Then $A=B_1^*\cup B_2^*\cup B_3^* \cup B_4^*$ and so $|B_1^*\cup B_2^*\cup B_3^* \cup B_4^*|\ge 5n-8$. Note that $B_1^*,  B_2^*, B_3^*, B_4^* $ are  pairwise disjoint.    
Suppose first that $B_1$ is red-complete to $B_3$ in $G$. By Lemma \ref{n,n+1 Lemma} applied to $B_3^*$ and $ B_1 \cup B_3$,    $|B_3^*| \le n-1$.  Thus  $|B_1^*| + |B_2^*| + |B_4^*| \ge 5n -8 - (n-1) = 4n - 7\ge 2n+5$ because $n\in\{6,7\}$.    
By symmetry, we may assume that $|B_1^*| +  |B_2^*| \ge n+3$.  We claim that $G[B_1^* \cup B_2^* \cup B_4^*]$ has no blue edges.  Suppose not. Let $uv$ be a blue edge in $G[B_1^* \cup B_2^* \cup B_4^*]$.   Since $|B_1^*| +  |B_2^*| \ge n+3$, let $x, y\in B_1^* \cup B_2^*$ be two distinct vertices that are different from $u$ and $v$.  If $u, v\in B_1^* \cup B_2^*$,  then we find a blue $C_{2n+1}$ with vertices 
\[
\begin{cases}
u, v, b_1, b_5, b_9, b_6, b_{10}, b_7, b_{11}, b_8, b_2, x, b_3, & \text{if $n=6$}\\
u, v, b_1, b_5, b_9, b_6, b_{10}, b_7, b_{11}, b_8, b_2, x, b_3, y, b_4, & \text{if $n=7$}
\end{cases}
\]
in order, a contradiction.   Thus we may assume that $v\in B_4^*$.   If $u \in B_1^*\cup B_2^*$,  then we find a blue $C_{2n+1}$ with vertices
\[
\begin{cases}
u, v, b_9, b_5, b_{10}, b_6, b_{11}, b_7, b_{12}, b_8, b_1, x, b_2, & \text{if $n=6$}\\
u, v, b_9, b_5, b_{10}, b_6, b_{11}, b_7, b_{12}, b_8, b_1, x, b_2, y, b_3, & \text{if $n=7$}
\end{cases}
\]
in order, a contradiction.  Thus  $u, v \in B_4^*$. But similarly, we obtain  a blue $C_{2n+1}$ with vertices 
\[
\begin{cases}
u, v, b_9, b_5, b_1, x, b_2, b_6, b_3, y, b_4, b_7, b_{10}, & \text{if $n=6$}\\
u, v, b_9, b_5, b_1, x, b_2, b_6, b_{10}, b_7, b_3, y, b_4, b_8, b_{11}, & \text{if $n=7$}
\end{cases}
\]
in order, a contradiction. This proves that  $G[B_1^* \cup B_2^* \cup B_4^*]$ contains no blue edges.  Since $|B_1^*| + |B_2^*| + |B_4^*| \ge 2n+5$ and $|A_p|\le n$, by Lemma \ref{Hamilton}, $G[B_1^* \cup B_2^* \cup B_4^*]$ has  a red $C_{2n+1}$, a contradiction.  Thus   $B_1$ must be  blue-complete to $B_3$.  Then $|B_1 \cup B_2 \cup B_3|\le 2n$, else we obtain a blue $C_{2n+1}$ in $G[B_1 \cup B_2 \cup B_3]$. By Lemma \ref{n,n+1 Lemma} applied to $B_2 \cup B_2^*$ and $B_1 \cup B_3$,  we see that  $|B_2^*| \le n-5$.  
If   $|B_1^*| \ge 3$, let   $x, y , z \in B_1^*$ be distinct vertices. Then we find a blue $C_{2n+1}$ with vertices  
\[
\begin{cases}
b_1, b_5, b_9, b_6, b_{10}, b_7, b_{11}, b_8, b_{12}, b_2, x, b_3, y, & \text{if $n=6$}\\
b_1, b_5, b_9, b_6, b_{10}, b_7, b_{11}, b_8, b_{12}, b_2, x, b_3, y, b_4, z, & \text{if $n=7$}
\end{cases}
\]
in order, a contradiction. Thus $|B_1^*| \le  2$. Similarly,   $|B_4^*| \le 2$. Therefore, 
\[
|B_3^*| = |G| - |X|-|B_1 \cup B_2 \cup B_3| - |B_1^* \cup B_2^* \cup B_4^*|\ge n\cdot 2^k+1-3k -2n -(n-5+2+2)   \ge 5n-7.
\]  
By Lemma \ref{n,n+1 Lemma} applied to $B_3^*$ and  $B_1 \cup B_3$, $G[B_3^*]$ contains no red edges.  But then by Lemma \ref{Hamilton} and the fact that $|A_p|\le n$ and $|B_3^*| \ge 5n-7$,  $G[B_3^*]$ must contain  a blue $C_{2n+1}$, a contradiction.  This proves that  if $|A_{p}| \le n$, then $|A_{p-2}| \le 3$.\qed\medskip

By Claim \ref{Ap4},   $|A_p|\ge4$ and so   $p-s\ge1$. Let 
\[
\begin{split}
B &:=  \{a_i \in \{a_1, \ldots, a_{p-1}\} \mid a_ia_p \text{ is colored blue in } \mc{R} \}\\
R &:= \{a_j \in \{a_1, \ldots, a_{p-1}\} \mid a_ja_p \text{ is colored red in } \mc{R} \}
\end{split}
\]
Then $|B|+|R|=p-1$. 
Let $B_G:= \bigcup_{a_i \in B} A_i$ and $R_G:=\bigcup_{a_j\in R} A_j$.

\begin{claim}\label{X'''}
If every vertex in $X$ is   neither $i$- nor $j$-complete to $V(G) \less X$ for two distinct colors $i, j\in[k]$, then $X'''= \emptyset$.
\end{claim}
\pf
Suppose $X''' \neq \emptyset$. We may assume that every vertex in $X$ is   neither red- nor blue-complete to $V(G) \less X$.   Then there exists at least one color $\ell \in [k]$ other than  red and blue such that  $|X^*_\ell|=3$. We claim that  $k\ge4$. Suppose $k=3$. Then $|G|=8n+1$.  We may assume    the third color is green. Then $|X|=|X^*_g|=3$.   By Claim~\ref{l-vertex} applied to  $X^*_{g}$  and $V(G) \less   X^*_{g}$, $|V(G) \less X^*_{g}|\le 4(2n-1) +(n-7)$. But then 
\begin{align*}
|G|  = |X^*_{g}  |+|V(G) \less   X^*_{g} |  \le  3+4(2n-1) +(n-7) <8n+1,   
\end{align*}
because $n\in\{6,7\}$, a contradiction. Thus  $k\ge4$, as claimed.   When choosing $X_1, X_2, \ldots, X_m$, let $q\in[m]$ be the smallest index such that for some color $\ell' \in [k]$ other than red and blue, $|X_{\ell'}^*\cap (X_1\cup  \dots \cup X_{q})|=3$. By the choice of $q$,  $|X_{j}^*\cap (X_1\cup  \dots \cup X_{q-1})|\le2$ for all $j\in[k]$. By the property (\ref{property b}) when choosing $X_1, X_2, \ldots, X_m$,  there are possibly two colors $q_1:=\ell'$ and  $q_2 \in[k]$ such that  $|X^*_{q_1}|=3$ and $|X^*_{q_2}|\le3$.  Since no vertex in $X$ is red- or blue-complete to $V(G) \less X$, we see that $|(X_1\cup\dots \cup X_{q}) \less (X^*_{q_1}\cup X^*_{q_2})|\le 2(k-4)$. By Claim~\ref{l-vertex} applied to  $X^*_{q_1}$  and $V(G) \less (X_1\cup  \dots \cup X_{q})$,  $|V(G) \less (X_1\cup \dots \cup X_{q})|\le (2n-1) \cdot 2^{k-1}+(n-7)$. But then 
\begin{align*}
|G| &=|(X_1\cup  \dots \cup X_{q}) \less (X^*_{q_1} \cup X^*_{q_2})|+|X^*_{q_1} \cup X^*_{q_2}|+|V(G) \less (X_1\cup \dots \cup X_{q})|\\
& \le  
2(k-4)+6+  [(2n-1) \cdot 2^{k-1}+(n-7)]\\
& < n \cdot 2^k+1
\end{align*}
for all $k \ge 4$, a contradiction.  \qed\\

\begin{claim}\label{B}
If $|A_p|\ge n$ and $|B|\ge 3$ (resp. $|R|\ge 3$), then $|B_G|\le 2n$ (resp. $|R_G|\le 2n$). 
\end{claim}

\pf
Suppose $|A_p|\ge n$ and $|B|\ge 3$ but $|B_G|\ge 2n+1$. By Claim \ref{n,n+1 Lemma},  $G[B_G]$ has no blue edges and no vertex  in $ X$ is blue-complete to $V(G)\less X$.  Thus all the edges of $\mathcal{R}[B]$ are colored red in $\mathcal {R}$. Let $q:=|B| $   and let $B:=\{a_{i_1}, a_{i_2}, \ldots, a_{i_q}\}$ with $|A_{i_1}|\ge |A_{i_2}|\ge\cdots\ge |A_{i_q}|$. Then $G[B_G]\less \bigcup_{j=1}^{q}E(G[A_{i_j}])$ is a complete multipartite graph with at least three parts.  If $|A_{i_1}|\le n$, then by Lemma~\ref{Hamilton} applied to $G[B_G]\less\bigcup_{j=1}^{q}E(G[A_{i_j}])$, $G[B_G]$ has a red $C_{2n+1}$, a contradiction.
 Thus  $|A_{i_1}|\ge n+1$.     Let $Q_b:=\{v\in R_G: v \text{ is blue-complete to } A_{i_1}\}$, and $Q_r:=\{v\in R_G: v \text{ is red-complete to } A_{i_1}\}$. Then $Q_b\cup Q_r=R_G$. 
 Let $Q: =(B_G\less A_{i_1})\cup Q_r \cup X^*_r$. Then $Q$ is red-complete to $A_{i_1}$ and $G[Q]$ must contain red edges,  because   $|B|\ge 3$ and all the edges of $\mathcal{R}[B]$ are colored red. By Claim \ref{n,n+1 Lemma} applied to $A_{i_1}$ and $Q$, $|Q|\le n$.
 Note that  $|A_p\cup Q_b|\ge|A_p|\ge|A_{i_1}|\ge n+1$ and $A_p\cup Q_b$ is blue-complete to $A_{i_1}$. By Claim \ref{n,n+1 Lemma} applied to $A_{i_1}$ and $A_p\cup Q_b$,   $G[A_p\cup Q_b]$ has   no  blue edges. Since no vertex  in $ X$ is blue-complete to $V(G)\less X$, we see that    $G[A_p\cup Q_b\cup (X'\less X^*_r)]$   has   no blue edges. By minimality of $k$,   $|A_p\cup Q_b\cup (X'\less X^*_r)|\le n\cdot 2^{k-1}$.
Suppose first that  $Q_r \cup X^*_r=\emptyset$.
 Then $Q_b=R_G$  and   $G[B_G\cup X'']$   has  no blue edges. By minimality of $k$, $|B_G\cup X''|\le n\cdot 2^{k-1}$. Since no vertex in $X$ is red- or blue-complete to $V(G) \less X$, by Claim~\ref{X'''}, $X'''=\emptyset$. But then
\[ |G| =|B_G\cup X''|+|A_p\cup Q_b\cup X'|  \le  n\cdot 2^{k-1}+n\cdot 2^{k-1} <n\cdot 2^{k}+1,\]
a contradiction. Thus $Q_r \cup X^*_r\ne\emptyset$. Since $|B|\ge3$, we see that $|B_G\less A_{i_1}|\ge2$. Thus  $n\ge |Q|\ge3$. 
 
 We next claim that either $|Q|\ge 4$ or $k\ge6$. Suppose $|Q|=3$ and $k\le 5$. Then $|Q_r \cup X^*_r|=1$ and $|B_G\less A_{i_1}|=2$.
 Suppose $k=3$.  We may assume that the third color is green.   Since $Q$ is red-complete to $A_{i_1}$, we see that $G[A_{i_1}]$ has neither red $C_{2n-2}$ nor a green $C_{2n+1}$.   By Theorem~\ref{R(Cm, Cl)}, $|A_{i_1}|\le R(C_{2n-2}, C_{2n+1})-1=4n-6$. But then
 \begin{align*}
|G| =|Q|+|A_{i_1}|+|A_p\cup Q_b| +|X^*_g|  \le 3+ (4n-6) +n\cdot 2^{3-1}+3   =8n<8n+1,
\end{align*}
a contradiction. Thus $k\in\{4,5\}$.  Then $|X'\less X_r^*|\le k-3$, else,   by Theorem~\ref{GR-even}, $|A_{i_1}|\le GR_{k-1}(C_{2n})-1\le (k-1)(n-1)+3n-1$. But then
 \begin{align*}
|G|&=|Q|+|A_{i_1}|+|A_p\cup Q_b\cup(X'\less X_r^*)| +|(X''\cup X''')\less X_r^*| \\
&\le 3+ [(k-1)(n-1)+3n-1] +n\cdot 2^{k-1}+2(k-2)  \\
& <n\cdot 2^k+1
\end{align*}
for all $k\in\{4,5\}$ and $n\in\{6,7\}$, a contradiction.  Thus $|X'\less X_r^*|\le k-3$, and so $|X''\less X^*_r|\le k-3$.  In particular, by Claim~\ref{X'''}, this implies $X'''=\es$.   By Claim \ref{l-vertex}  applied to $Q$ and  $ A_{i_1}$,  $|A_{i_1}|\le (2n-1) \cdot 2^{k-2}+(n-7)$. But then 
 \begin{align*}
|G|&=|Q|+|A_{i_1}|+|A_p\cup Q_b\cup (X'\less X^*_r)| +|X''\less X_r^*|\\
  &\le  3+ [(2n-1) \cdot 2^{k-2}+(n-7)]+n\cdot 2^{k-1}+(k-3)\\
  & <n\cdot 2^k+1
 \end{align*}
 for $k \in \{4, 5\}$, a contradiction.  This proves that  either $|Q|\ge 4$ or $k\ge6$, as claimed. 
 
 Note that $G[A_{i_1}]$ has no blue edges and $|(X'' \cup X''')\less X^*_r|\le 2(k-2)$.
 By Claim \ref{l-vertex}  applied to $Q$ and  $ A_{i_1}$, we see that
  \[
   |A_{i_1}|\le\begin{cases} 
  n \cdot 2^{k-2}  &  \text{if } |Q| = n\\
n \cdot 2^{k-2} +2  &  \text{if } |Q| = n-1\\
  (21-2n) \cdot 2^{k-2} + (5n - 31)  &  \text{if } |Q| = n-2\\
   11 \cdot 2^{k-2} + (n-7)  &  \text{if } |Q| = n-3\\
13 \cdot 2^{k-2}  &  \text{if } |Q| = n-4 \text{ and } n=7.
\end{cases}
\]
But then
\begin{align*}
|G|&=|Q|+|A_{i_1}|+|A_p\cup Q_b\cup (X'\less X^*_r)| +|(X'' \cup X''')\less X^*_r|\\
  &\le\begin{cases} 
  n+n \cdot 2^{k-2}+n\cdot 2^{k-1}+2(k-2)  &  \text{if } |Q| = n\\
  (n-1)+(n \cdot 2^{k-2} +2) +n\cdot 2^{k-1}+2(k-2)  &  \text{if } |Q| = n-1\\
  (n-2)+[(21-2n) \cdot 2^{k-2} +(5n-31)] +n\cdot 2^{k-1}+2(k-2)  &  \text{if } |Q| = n-2\\
  (n-3)+ [11 \cdot 2^{k-2}+(n-7)] +n\cdot 2^{k-1}+2(k-2)  &  \text{if $|Q| = n-3$ and $n=7$}\\
  3+ [(2n-1) \cdot 2^{k-2}+(n-7)] +n\cdot 2^{k-1}+2(k-2)  &  \text{if $|Q| = 3$ and $k\ge6$}.\\
\end{cases}
\end{align*}
In each case, we have $|G|< n\cdot 2^{k}+1$,  a contradiction.
This proves that  if  $|A_p|\ge n$ and $|B|\ge 3$, then $|B_G|\le 2n$. Similarly, one can prove that if $|A_p|\ge n$ and $|R|\ge 3$, then $|R_G|\le 2n$. \qed\\

\begin{claim}\label{p}
$p\le 2n+1$.
\end{claim}
\pf
Suppose $p \ge 2n+2$. Then $|B| + |R| = p-1 \ge 2n+1$. We claim that $|A_p| \le n$. Suppose $|A_p| \ge n+1$. We may assume that $|B| \ge |R|$.  Then $|B_G| \ge |B| \ge n+1$. By Claim~\ref{B},  $|B_G| \le 2n$, and by Claim~\ref{n,n+1 Lemma}, $G[A_p]$ has no blue edges and $X_b^* = \es$.  Then $|X''\cup X'''|\le 2(k-1)$. If  $|R_G|\ge n+1$, then by Claim~\ref{n,n+1 Lemma}, neither $G[R_G]$ nor  $G[A_p]$ has red edges and $X_r^* = \es$. By Claim~\ref{X'''}, $X'''=\emptyset$. 
Note that $G[A_p \cup X']$ has neither red nor blue edges, and $G[R_G \cup X'']$ has no red edges.  Then by minimality of $k$,
\[
|G| = |A_p \cup X'| + |B_G| + |R_G \cup X''| + |X'''| \le n \cdot 2^{k-2} + 2n + n \cdot 2^{k-1} < n \cdot 2^k + 1
\]
for all $k \ge 3$, a contradiction. Thus, $|R_G|\le n$. Then for all $k \ge 3$,
\[
|A_p\cup X'|=|G|-|B_G|-|R_G|-|X'' \cup X'''|\ge n\cdot 2^k+1- 2n -n-2(k-1)>n\cdot 2^{k-1}+1.
\]
Since $G[A_p\cup X']$ has no blue edges, by the choice of $k$, $G[A_p\cup X']$ has a monochromatic $C_{2n+1}$, a contradiction. This proves that $|A_p| \le n$, as claimed.

Note that by Claim \ref{Ap4}, $|A_p|\ge 4$. Additionally, Claims~\ref{Ap-2} and \ref{Ap-8} give $|A_{p-2}| \le 3$ and $|A_{p-8}| \le 2$ with $|A_{p-4n+12}| \le 1$, respectively.  Therefore, $k = 3$, $|G| = 8n+1$ and $|X| \le 9$.  Because $n \in \{6,7\}$, 
\[
|B_G| + |R_G| = |G| - |A_p| - |X| \ge (8n+1) - n-9 =7n-8 > 6n-3 \ge GR_3(C_{2n})
\]  
by Theorem~\ref{GR-even}.  Therefore, $|X| \le 6$, otherwise we find a monochromatic $C_{2n+1}$.  Recalculating the above inequality with this fact, we obtain
\[
|B_G| + |R_G| = |G| - |A_p| - |X| \ge (8n+1) - n-6 =7n-5.
\]
Thus at least one of $|B_G| \ge 3n + 1$  or  $|R_G| \ge 3n+1$, so we may assume $|B_G| \ge 3n+1$ in what follows, because the argument is identical if $|R_G| \ge 3n+1$.   We next prove that $4 \le |A_p|\le n$ is impossible.  

Suppose first that $|A_p| \ge 5$ and let $B^* \subseteq B_G$ be a minimal set such that $G[B_G \less B^*]$ has no blue edges.  If $|B^*| \le 2n-10$, then $|B_G \less B^*| \ge (3n+1) - (2n-10)= n+11 \ge 2n + 4$, because $n \in \{6,7\}$. Therefore, $|B \less B^*| \ge 3$, and all edges of $\mc{R}[B\less B^*]$ are colored red, so that by Lemma \ref{Hamilton} we find a monochromatic $C_{2n+1}$, a contradiction.  Thus, $|B^*| \ge 2n - 9$.  Define the family of graphs 
\begin{align*}
\mc{H}_1 &:= \{(2n-9)K_2, (2n-11)K_2 \cup P_3, (15-2n)K_2 \cup 2P_{2n-11}, \\
	& \hspace{1in} 2P_{n-5} \cup P_4, P_{2n-11} \cup P_4, (15-2n)K_2 \cup P_{4n-23}, P_{2n-8}\}
\end{align*}
It follows that $G[B_G]$ contains a blue $H \in \mc{H}_1$, so that along with the vertices in $A_p$, we find a blue $C_{2n+1}$, a contradiction.

Therefore suppose $|A_p| = 4$.  Then $|B_G| + |R_G| = |G| - |A_p| - |X| \ge 8n + 1 - 4 - 6\ge  8n -9$, so that $|B_G| \ge 4n-4$.  Let $B^*$ be defined as above.  If $|B^*| \le 2n-5$, then $|B_G \less B^*| \ge (4n-4) - (2n-5)= 2n + 1$, and thus $|B \less B^*|\ge 3$.  Since $\mc{R}[B \less B^*]$ contains only red edges, by Lemma \ref{Hamilton}, there is a red $C_{2n+1}$, a contradiction.  Thus, $|B^*| \ge 2n-4$.  Define the family of graphs
\begin{align*}
\mc{H}_2:=\{ (2n-4)K_2, (14-n)K_2 \cup P_{3n-17}, (20-2n)K_2 \cup 2P_{2n-11}, 8K_2 \cup P_{2n-11}\}.
\end{align*}
Let $M$ denote a matching of size $m \ge 0$.  For any $H \in \mc{H}_2$, let $H':=H \cup M$.  It follows that $G[B_G]$ contains a blue $H'$, where $m$ is chosen to be as large as possible.  Then removing at most two vertices, say $x, y \in V(H)$ from the longest blue subpaths in $H$, we obtain $M':=H' \less \{x, y\}$, which is a matching of size $m' \ge 6$.  Denote the edges in $M'$ by $u_iv_i$, for all $i \in [m']$.  Put another way, this means the blue edges in $G[B_G \less \{x, y\}]$ induce a blue matching.  Let us define a new Gallai partition of $G[B_G \less \{x, y\}]$ in the following manner.  If $|A_{i_j}| = |A_{i_ \ell}| = 1$ for some pair $j, \ell \in [q]$, and if $A_{i_j}$ is blue-complete to $A_{i_ \ell}$, then create the new part $A_{i_s}:= A_{i_j} \cup A_{i_\ell}$, so that $|A_{i_s}| = 2$, where $s \in [q']$ and $q' \le q$;  otherwise, define $A_{i_j}$ to be the same.  By construction, only red edges appear between any two parts of this modified partition. We may assume $A_{i_1}, \ldots, A_{i_t}$ are all parts of the modified Gallai partition of $B_G \less \{x, y\}$ containing blue edges.  Because $m'\ge 6$, we see that $\bigcup\limits_{j=1}^t |A_{i_j}| \ge 12$, and because $|A_p| = 4$, we also have $t \ge 3$.  In particular, if $\bigcup\limits_{j=1}^t |A_{i_j}| \ge 2n+1$, we are done by Lemma~\ref{Hamilton} because $G\left[\bigcup\limits_{j=1}^t A_{i_j} \right] - \bigcup\limits_{j=1}^t E(A_{i_j})$ is a complete multipartite graph containing only red edges.  Thus we may assume $12 \le \sum\limits_{j=1}^t |A_{i_j}| \le 2n$.  Note that $|B_G \less \{x, y\}| - \sum\limits_{j=1}^t |A_{i_j}| \ge (4n-4) - 2 - 2n = 2n-6$.  Define $r:=2n+1 - \sum\limits_{j=1}^t |A_{i_j}|$, and choose distinct vertices $v_1, \ldots, v_r \in B_G \less \left(\{x, y\} \cup \bigcup\limits_{j=1}^t A_{i_j}\right)$.  Because $v_1, \ldots, v_r \not\in \bigcup\limits_{j=1}^t A_{i_j}$, we see that $\{v_1, \ldots, v_r\}$ is red-complete to $\bigcup\limits_{j=1}^t A_{i_j}$, again yielding a red $C_{2n+1}$ by Lemma~\ref{Hamilton}, again forcing a contradiction.   \qed\\

\begin{claim}\label{Ap}
$|A_p|\ge n+1$.
\end{claim}

\pf  Suppose $|A_p|\le n$.  Then $p\ge 9$ because $|G|\ge8n+1$. 
By Claim \ref{p}, $9\le p \le 2n+1$.  We may assume that $a_pa_{p-1}$ is colored blue in $\mc{R}$. Then $|A_p\cup A_{p-1}\cup X_b^* |\le 2n$, else $|X_b^*|\ge1$ and so $G[A_p\cup A_{p-1}\cup X_b^*]$ has a blue $C_{2n+1}$, a contradiction.   It follows that  $|A_p\cup A_{p-1}\cup X|=|A_p\cup A_{p-1}\cup X_b^*|+|X\less X_b^*|\le 2n+3(k-1)$. By Claim \ref{Ap-2} and Claim \ref{Ap-8},   $|A_{p-2}|\le3$ and  $|A_{p-8}|\le 2$. 
But then  
\[
\begin{split}
|G| 	&=|A_p\cup A_{p-1}\cup X|+\sum_{i=p-7}^{p-2} |A_i|+ \sum_{i=1}^{p-8} |A_i| \\
& \le  [2n+3(k-1)]+18 +2(2n+1-8)\\
&= 6n+3k+1\\
& <n\cdot 2^k+1,   
\end{split}
\]
for    $n\in\{6,7\}$ and all $k\ge3$, a contradiction. \qed \\
 
\begin{claim}\label{A_{p-2}}
$|A_{p-2}|\le n$.
\end{claim}

\pf Suppose $|A_{p-2}|\ge n+1$. Then $n+1 \le |A_{p-2}| \le |A_{p-1}| \le |A_p|$ and so  $\mc{R}[\{a_{p-2}, a_{p-1},a_p\}]$ is not a monochromatic triangle in $\mc{R}$ (else $G[A_p\cup A_{p-1}\cup A_{p-2}]$ has a   a monochromatic $C_{2n+1}$).   Let $B_1$, $B_2$, $B_3$ be a permutation of $A_{p-2}$, $A_{p-1}$, $A_p$ such that $B_2$ is, say blue-complete,  to $B_1 \cup B_3$ in $G$. Then $B_1$ must be  red-complete to $B_3$ in $G$.  By Claim~\ref{n,n+1 Lemma}, $X_r^* = \es$ and $X_b^* = \es$.  By Claim~\ref{X'''}, $X'''=\emptyset$.   Let $A:=V(G)\less (B_1\cup B_2\cup B_3\cup X'\cup X'')$.  By Claim~\ref{n,n+1 Lemma} again, $G[B_2]$ has no blue edges, and neither $G[B_1\cup X']$ nor $G[B_3\cup X'']$ has red or blue edges.  By minimality of $k$,  $|B_1\cup X'|\le n\cdot 2^{k-2}$ and $|B_3\cup X''|\le n\cdot 2^{k-2}$. It follows that 
 $|A\cup B_2|=|G|-|B_1\cup X'|-| B_3\cup  X''|\ge n\cdot 2^{k-1}+1$.  By minimality of $k$, $G[A\cup B_2]$ must have blue edges. By Claim~\ref{n,n+1 Lemma},  no vertex in $A$ is red-complete to $B_1\cup B_3$ in $G$,  and no vertex in $A$ is blue-complete to $B_1\cup B_2$ or $B_2\cup B_3$ in $G$. This implies that 
$A$ must be  red-complete to $B_2$ in $G$. It follows that $G[A]$ must contain a blue edge, say $uv$.    
 Let $b_1, \ldots, b_{n-1}\in B_1$, $ b_{n}, \ldots, b_{2n-2}\in B_2$, and $b_{2n-1}\in B_3$.   If $\{u,v\}$   is blue-complete to $ B_1$, then we obtain a blue $C_{2n+1}$ with vertices $b_1, u, v, b_2, b_n, b_{2n-1}, b_{n+1}, b_3, b_{n+2}, \ldots, b_{n-1}, b_{2n-2}$
   in order, a contradiction.  Thus $\{u,v\}$   is not blue-complete to $ B_1$. Similarly, $\{u,v\}$   is not blue-complete to $ B_3$.
Since no vertex in $A$ is red-complete to $B_1\cup B_3$,  we may assume that $u $ is blue-complete to $ B_1$ and $v $ is blue-complete to $ B_3$. But then we obtain a blue $C_{2n+1}$ with vertices $b_1, u, v, b_{2n-1}, b_n, b_2, b_{n+1},\ldots,  b_{n-1}, b_{2n-2}$ in order.  \qed \\

For the remainder of the proof, let  $B_G^*:=B_G\cup X_b^*$   and   $R_G^*:=R_G\cup X_r^*$.

\begin{claim}\label{BR*_G}
$|B_G|\ge4$ or $|R_G|\ge4$. 
\end{claim}

\pf Suppose $|B_G|\le3$ and $|R_G|\le3$.  Since $p\ge2$, we see that  $B_G\ne\emptyset$ or $R_G\ne\emptyset$.   By maximality of $m$ (see condition (\ref{property a}) when choosing $X_1, X_2, \ldots, X_m$), $B_G\ne \emptyset$, $R_G\ne \emptyset$, and $B_G$ is neither red- nor blue-complete to $R_G$ in $G$. But then, since $|B_G|\le3$ and $|R_G|\le3$, by maximality of $m$ again (see condition (\ref{property b}) when choosing $X_1, X_2, \ldots, X_m$), $B_G =\emptyset$ and   $R_G=\emptyset$,  a contradiction.\qed\medskip

\begin{claim}\label{one}
$2\le p-s\le 8$.
\end{claim}

\pf  By Claim \ref{Ap-8}, $|A_{p-8}|\le2$ and so     $p-s\le 8$. Suppose $p-s\le1$. Then $p-s=1$ because $p-s\ge1$.  Thus $|A_i|\le2$ for all $i\in[p-1]$ by the choice of $p$ and $s$.  By Claim \ref{p}, $p \le 2n+1$.    Then $|B_G \cup R_G|\le 2(p-1)$ and so $|B^*_G \cup R^*_G| \le 2(p-1)+3+3=2(p+2)\le 4n+6$. We may assume that $|B^*_G| \ge |R^*_G|$.  If $|R^*_G|\ge n$, then $|B^*_G|\ge n$. By Claim \ref{Ap} and Claim \ref{n,n+1 Lemma}, $G[A_p]$ has neither  blue nor red edges. By minimality of $k$, $|A_p  | \le n \cdot 2^{k-2}$. But then 
\[
|G|=|B^*_G \cup R^*_G|+|A_p |+ |X \less (B^*_G \cup R_G^*)|\le (4n+6) + n \cdot 2^{k-2} +3(k-2)< n\cdot 2^k+1
\]
for all $k\ge3$, a contradiction.  Thus $|R^*_G|\le n-1$. We claim that   $|B^*_G|\le 2n+3$. This is trivially true if $|B|\le n$.   If $|B|\ge n+1$, then  $|B_G|\le 2n$  by Claim~\ref{B}.  Thus $|B^*_G|\le 2n+3$, as claimed.   If $|B_G^*| \ge n-1$, then applying Claim \ref{l-vertex}(i,ii) to   $B_G^*$ and $A_p$ implies that 
 \[
  |B_G^*| + |A_p|\le\begin{cases}  
 (n-1)+(n \cdot 2^{k-1} +2), &  \text{if } |B_G^*| = n-1\\
  (2n+3)+n \cdot 2^{k-1}, &  \text{if } |B_G^*| \ge n.
 \end{cases}
 \] 
In either case, $ |B_G^*| + |A_p|\le (2n + 3)+n \cdot 2^{k-1}$. But then 
\[
|G| = |R^*_G| + |B_G^*| + |A_p| + |X \less (B^*_G \cup R_G^*)| \le  (n-1) + [(2n+ 3)+n \cdot 2^{k-1}] +  3(k-2) <n \cdot 2^k + 1,
\]
for all $k\ge3$  and $n\in\{6,7\}$, a contradiction. 
Thus $|R_G^*| \le |B_G^*| \le n-2$.  If $|B_G^*| = n-2$, then by Claim \ref{l-vertex}(iii), $|A_p| \le (21-2n) \cdot 2^{k-1} + (5n-31)$.  But then 
\[
|G| = |R^*_G| + |B_G^*| + |A_p| + |X \less (B^*_G \cup R_G^*)| \le  2(n-2) + [(21-2n) \cdot 2^{k-1} + (5n-31)] +  3(k-2) <n \cdot 2^k + 1,
\]
for all $k \ge 3$ and $n \in \{6,7\}$, a contradiction.   Thus $|R_G^*| \le |B_G^*|\le  n-3$. 
By Claim \ref{BR*_G},  $|R_G^*| \le |B_G^*| =|B_G|=4$  and  $n=7$.   By Claim \ref{l-vertex}(iv), $|A_p| \le 11\cdot 2^{k-1}$.  But then 
\[
|G| = |R^*_G| + |B_G^*| + |A_p| + |X \less (B^*_G \cup R_G^*)| \le  4 + 4 + 11\cdot 2^{k-1}+  3(k-2) <7 \cdot 2^k + 1,
\]
for all $k \ge 3$, a contradiction.\qed \\

  By Claim~\ref{one},   $2 \le p-s \le  8$ and so $|A_{p-1}|\ge3$.  We may now assume that $a_pa_{p-1}$ is colored blue in $\mc{R}$. Then   $a_{p-1}\in B$ and so $A_{p-1}\subseteq B_G$. Thus $|B^*_G|\ge|B_G|\ge|A_{p-1}|\ge3$.

\begin{claim}\label{R*_G}
$|R^*_G|\le 2n$.
\end{claim}

\pf 
Suppose   $|R^*_G| \ge 2n+1$.  By Claim~\ref{Ap}, $|A_p| \ge n+1$. By Claim~\ref{n,n+1 Lemma},  $G[R^*_G]$ has no red edges.  Thus $|R^*_G|=|R_G|$ and so $X_r^* = \es$.  In particular, all the edges in $\mc{R}[R]$ are colored blue. By  Claim~\ref{B}, $|R| \le 2$.  By Claim~\ref{A_{p-2}},  $|A_{p-2}| \le n$.    Since $A_{p-1}\cap R_G=\emptyset$ and $|R_G| \ge 2n+1$, we see that $|R| \ge 3$, a contradiction.  \qed\\

\begin{claim}\label{A_{p-1}}
$|A_{p-1}|\le n$.
\end{claim}
 
Suppose   $|A_{p-1}| \ge n+1$.  
Then $|B_G| \ge |A_{p-1}|\ge n+1$.  By Claim~\ref{n,n+1 Lemma},  neither $G[A_p]$ nor $G[B_G]$ has blue edges, and $X_b^* = \es$.  Thus $|X|\le 3(k-1)$.  We   claim that $X''' = \es$.  Suppose $X''' \ne \es$.  By Claim \ref{X'''}, $|X^*_i|\ge1$ for every  color $i\in[k]$ other than blue,  and   $|X_j^*|=3$ for some  color $j\in[k]$ other than blue.  Then by Claim \ref{l-vertex}(iv, v) applied to $X_j^*$ and $V(G) \less X$, $|V(G) \less X| \le (2n-1) \cdot 2^{k-1} + n-7$.  Thus $|X|\ge 3k-4$, else,   
\[
|G| = |V(G) \less X| + |X| \le [(2n-1) \cdot 2^{k-1} + n-7]+ 3k-5 < n \cdot 2^k + 1
\]
for all $k \ge 3$, a contradiction. 
     We claim that  $k\ge4$. Suppose $k=3$.  We may assume that the third color is green. Since $|X|\ge 3k-4=5$, we have $|X_r^*|\ge2$ and $|X_g^*|\ge2$.  By Claim \ref{n,n+1 Lemma} applied to $A_p$ and $R^*_G$, $|R^*_G|\le n$. Thus $|A_p|+|B_G|=|G|-|R^*_G|-|X_g^*|\ge 8n+1-n-3=7n-2$. Thus either $|A_p|\ge 3n+2$ or   $|B_G|\ge 3n+2$. We may assume that $|A_p|\ge 3n+2$. By Theorem \ref{even cycles}, $G[A_p]$ has either a red or a green $C_{2n}$. Thus either $G[A_p\cup X_r^*]$ has a red $C_{2n+1}$ or $G[A_p\cup X_g^* ]$ has a green $C_{2n+1}$, a contradiction. 
  Thus  $k\ge4$, as claimed.  
Since $|X|\ge 3k-4$,  by Claim \ref{X*}, we may assume that $2\le |X^*_g|\le3$, and  $|X^*_i|=3$ for every  color $i\in[k]$ other than blue and green. 
When choosing $X_1, X_2, \ldots, X_m$, let $q\in[m]$ be the smallest index such that for some color $\ell \in [k]$ other than blue, $|X_{\ell}^*\cap (X_1\cup  \dots \cup X_{q})|=3$. By the choice of $q$,  $|X_{j}^*\cap (X_1\cup  \dots \cup X_{q-1})|\le2$ for all $j\in[k]$. By the property (\ref{property b}) when choosing $X_1, X_2, \ldots, X_m$,  there are possibly two colors $q_1, q_2 \in[k]$ such that  $q_1=\ell$,  $|X^*_{q_1}\cap (X_1\cup  \dots \cup X_{q})|=3$ and $|X^*_{q_2}\cap (X_1\cup  \dots \cup X_{q})|\le3$.  Since $X_b^* = \es$, $k\ge4$  and $|X^*_i|=3$ for every  color $i\in[k]$ other than blue and green, we see that $q<m$ and so $|(X_1\cup\dots \cup X_{q}) \less (X^*_{q_1}\cup X^*_{q_2})|\le 2(k-4)$. By Claim~\ref{l-vertex} applied to  $X^*_{q_1}$  and $V(G) \less (X_1\cup  \dots \cup X_{q})$, 
$|V(G) \less (X_1\cup \dots \cup X_{q})|\le (2n-1) \cdot 2^{k-1}+(n-7)$. But then 
\begin{align*}
|G| &=|(X_1\cup  \dots \cup X_{q}) \less (X^*_{q_1} \cup X^*_{q_2})|+|X^*_{q_1} \cup X^*_{q_2}|+|V(G) \less (X_1\cup \dots \cup X_{q})|\\
& \le  
2(k-4)+6+ [(2n-1) \cdot 2^{k-1}+(n-7)]
 \\
& < n \cdot 2^k+1
\end{align*}
for all $k \ge 4$, a contradiction. 
This proves that  $X''' = \es$,  as claimed. Thus  $|X| \le 2(k-1)$.  

  Since neither $G[A_p]$ nor $G[B_G]$ has blue edges and $X_b^* = \es$, we see that neither $G[A_p\cup X']$ nor $G[B_G\cup X'']$ has blue edges. By the choice of $k$, $|A_p\cup X'|\le n \cdot 2^{k-1}$ and $|B_G\cup X''|\le n \cdot 2^{k-1}$. We claim that $G[R_G]$ has blue edges. 
  
  Suppose  $G[R_G]$ has no blue edges. Then $G[A_p\cup R_G\cup X']$ has no blue edges. By the choice of $k$, $|A_p\cup R_G\cup X'|\le n \cdot 2^{k-1}$. But then $|B_G\cup X''|=|G|-|A_p\cup R_G\cup X'|\ge n \cdot 2^{k-1}+1$, a contradiction. Thus $G[R_G]$ has blue edges, as claimed. Then  $|R_G|\ge2$.  By Claim~\ref{R*_G},   $2\le |R_G|\le |R^*_G|\le 2n$.   Suppose   $|R^*_G| \ge n-1$. We claim that $|A_p\cup (X'\less X^*_r)|+|R^*_G| \le n \cdot 2^{k-2}+\max\{2n, k+n-1\}$. If $|R^*_G| \ge n$, then  by Claim~\ref{n,n+1 Lemma}, $G[A_p]$ has no red edges and so $G[A_p\cup (X'\less X^*_r)]$ has no red edges.    By the choice of $k$, $|A_p\cup (X'\less X^*_r)|\le n \cdot 2^{k-2}$ and so $|A_p\cup (X'\less X^*_r)|+|R^*_G|\le n \cdot 2^{k-2}+2n$. If  $  |R_G^*| = n-1$, then applying Claim \ref{l-vertex}(ii) to   $R^*_G$ and $A_p$, $|A_p| \le n \cdot 2^{k-2}+2$. Thus 
$|A_p\cup (X'\less X^*_r)|+|R^*_G|\le (n \cdot 2^{k-2}+2)+(k-2) +(n-1)=n \cdot 2^{k-2} +k+n-1$.  Thus $|A_p\cup (X'\less X^*_r)|+|R^*_G|\le n \cdot 2^{k-2}+\max\{2n, k+n-1\}$, as claimed. But then
\[  
 |G|=|A_p\cup (X'\less X^*_r)|+|R^*_G|+|B_G\cup (X''\less X^*_r)| \le (n\cdot 2^{k-2}+\max\{2n, k+n-1\})+ n\cdot 2^{k-1}   < n \cdot 2^k+1,
\]
for all $k\ge3$, a contradiction. 

Next, suppose $|R_G^*| = n-2$.  Then by applying Claim \ref{l-vertex}(iii) to $R_G^*$ and $A_p$, $|A_p| \le (21-2n)\cdot 2^{k-2} + (5n-31)$.  But then

 \[ 
 \begin{split} 
 |G|	&\le |A_p|+|B_G\cup X''|+|R^*_G|+|X'\less X^*_r|\\
 		&\le [(21-2n)\cdot 2^{k-2} + (5n-31)]+ n\cdot 2^{k-1} +(n-2)+ (k-2)\\
 		&< n \cdot 2^k+1,
 \end{split}
 \]
for all $k \ge 3$, a contradiction.  Thus $|R_G^*| \le n-3$. If $|R_G^*| = 4$, then  $n=7$ and so  by Claim \ref{l-vertex}(iv) applied to $R_G^*$ and $A_p$, $|A_p|\le 11\cdot 2^{k-2}$. But then 

 \[ 
 \begin{split} 
 |G|	&\le |A_p|+|B_G\cup X''|+|R^*_G|+|X'\less R^*_G|\\
 		&\le 11\cdot 2^{k-2}+ 7\cdot 2^{k-1} +4 + (k-2)\\
 		&< 7 \cdot 2^k+1,
 \end{split}
 \]
for all $k \ge 3$, a contradiction. Therefore, $|R_G^*| \le 3$.

Let $xy$ be a blue edge in $G[R_G]$. This is possible because $G[R_G]$ has blue edges.  We claim that either $x$ or $y$ is red-complete to $B_G$. Suppose there exist $x', y'\in B_G $   such that $xx'$ and $yy'$ are colored blue. Then $x'=y'$, else we obtain a blue $C_{2n+1}$ by Claim~\ref{n,n+1 Lemma} applied to $B_G$ and $A_p \cup \{x, y\}$. Thus  $x'$ is the  unique vertex in $  B_G$ such that    $\{x, y\}$  is  red-complete to $B_G\less x'$ in $G$ and $xx'$, $ yx'$ are colored blue.  Then there exists $i\in[s]$ such that $A_i=\{x'\}$.   Since $G[B_G]$ has no blue edges, we see that 
 $\{x, y, x'\}$ must be red-complete to $B_G\less x'$ in $G$.  

Now, if $|R_G^*| = 3$, let $R_G^* = \{x, y, z\}$.  If either $zx$ or $zy$ is blue, then  $X_r^* = \es$ and by the above reasoning, $z$ is also red-complete to $B_G \less x'$.  The same is true if $z \in X_r^*$.  By Claim~\ref{l-vertex}(iii, iv), $|B_G \less x'| \le (2n-3) \cdot 2^{k-2} + (n-7)$.  Again, by Claim~\ref{l-vertex}(iv, v), $|A_p| \le (2n-1) \cdot 2^{k-2} + (n-7)$.  But then
 \begin{align*}
 |G| &= |A_p| + |B_G \less x'| + |R_G^* \cup x'| + |X|\\
 		&\le [(2n-1) \cdot 2^{k-2} + (n-7)] + [(2n-3) \cdot 2^{k-2} + (n-7)] + 4 + 2(k-2)\\
 		&<n \cdot 2^k + 1
 \end{align*}
 for all $k \ge 3$, a contradiction.  Therefore, we may assume both $zx$ and $zy$ are red, but that $z \not\in X_r^*$.
 
In what follows, we now assume $2 \le |R_G| \le |R_G^*| \le 3$. By Claim~\ref{l-vertex}(iv, v) applied to  $\{x, y, x'\}$ and  $B_G\less x'$,    $|B_G\less x'|\le (2n - 1) \cdot 2^{k-2} +n-7$.   Note that $G[A_p \cup X' \cup \{x, z\}]$ has no blue edges if $|R_G^*| = 3$, and similarly $G[A_p \cup X' \cup  \{x\}]$ if $|R_G| = 2$.   Then $|X''|\ge k-2$, else,
   \[  
 |G|= |A_p\cup X' \cup \{x, z\}|+|B_G\less x'|+|\{y, x'\}|+|X''|\le  n\cdot 2^{k-1} +[(2n-1) \cdot 2^{k-2} + n-7]+2+(k-3)< n \cdot 2^k+1,
  \]
 for all $k \ge 3$, a contradiction. Since $2\le |R_G|\le |R^*_G|\le3$, we see that  $|X_r^*|\le1$. 
  It follows that $|X_i^*|=2$ for all colors $i\in[k]$ other than red and blue. Then neither $G[A_p]$ nor  $G[B_G\less \{x'\}]$ has  a monochromatic  $C_{2n-1}$ in any color $i\in [k]$ other than red and blue. Clearly, neither $G[A_p]$ nor  $G[B_G\less\{x'\}]$ has  red $C_{2n-1}$ because $\{x, y\}$ is red-complete to both  $A_p$ and $B_G\less\{x'\}$. By Theorem~\ref{C9C11} for $n=6$ and Theorem~\ref{C13C15} for $n=7$ (because although proved simultaneously here, the proof for $n=6$ is independent), $|B_G\less x'| \le (n-1) \cdot 2^{k-1}$ and $|A_p| \le (n-1) \cdot 2^{k-1}$.   But then 
 \[  
 |G|= |A_p|+|B_G\less x'|+|R_G^*\cup\{  x'\}|+|X \less X_r^*| \le  (n-1)\cdot 2^{k-1} + (n-1)\cdot 2^{k-1}+4+2(k-2)< n \cdot 2^k+1,
  \]
for all $k \ge 3$, a contradiction.  This proves that  either $x$ or $y$ is red-complete to $B_G$.  We may assume that $x$ is red-complete to $B_G$.

Suppose $|R_G|=2$. Then $R_G=\{x,y\}$ and $|X_r^*|\le 1$.    It follows that neither $G[A_p\cup\{y\}\cup X']$ nor $G[B_G\cup\{x\}\cup X'']$ has blue edges. By minimality of $k$, $|A_p\cup \{y\}\cup X'|\le n \cdot 2^{k-1}$ and $|B_G\cup\{x\}\cup X''|\le n \cdot 2^{k-1}$. But then  $|G|=|A_p\cup \{y\}\cup X'|+|B_G\cup\{x\}\cup X''|\le n\cdot 2^{k-1}+ n\cdot 2^{k-1}< n\cdot 2^k+1$ for all $k\ge3$, a contradiction. 
Thus  $|R_G| = |R_G^*| = 3$. Then $X_r^*=\es$ and   $G[A_p]$  has no     red $C_{2n}$.  
Clearly,  $|X'| \le  k-2$.  We claim that $|X'| \le  k-3$.  Suppose $|X'| = k-2$. Then  $|X_i^*|\ge1$ for all color $i\in[k]$ other than red and blue. Thus     $G[A_p]$    has   no  monochromatic  $C_{2n}$ in any colors $i\in [k]$ other than blue. 
 Since   $G[A_p]$    has  no  blue edges, by Theorem \ref{GR-even}, $|A_p| \le (n-1)(k-1) + 3n-1$.  Then $k=3$, else, 
\[
|G| = |A_p|+|B_G\cup X''|+|R_G|+|X' |\le [(n-1)(k-1) + 3n-1] + n\cdot 2^{k-1} + 3 + (k-2) < n \cdot 2^k+1
\]
for all $k \ge 4$.    Since $R_2(C_{2n}) = 3n-1$, we see that   $|A_p| \le 3n-2$ if $k=3$.  But then 
\[
|G| = |A_p| + |B_G \cup X''| + |R_G| + |X'| \le (3n-2) + 4n + 3 + 1 = 7n+2 < 8n+1.
\]
Thus $|X''|\le |X'| \le k-3$, as claimed.   Since $x$ is red-complete to $B_G$, we see that $G[B_G\cup \{x\}\cup X'']$ has no blue edges. By minimality of $k$, $|B_G\cup\{x\}\cup X''|\le  n\cdot 2^{k-1}$.   By Claim \ref{l-vertex}  applied to $R_G$ and $A_p$,  $|A_p|\le (2n-1) \cdot 2^{k-2} + n-7$. But then   
\[
|G| =  |A_p|+|B_G\cup\{x\}\cup X''|+|R_G\less x|+|X'|\le [(2n-1) \cdot 2^{k-2} + n-7] + n\cdot 2^{k-1} + 2 + (k-3) < n \cdot 2^k+1
\]
 for all $k \ge 3$, a contradiction.  Hence, $|A_{p-1}|\le n$.
 \qed\\

By Claim \ref{A_{p-1}},  $|A_{p-2}|\le |A_{p-1}| \le n$. Then  $|B_G| \le 2n$, because this is trivially true when $|B|\le 2$, and follows from Claim \ref{B} when $|B| \ge 3$.  By Claim~\ref{R*_G}, $ |R_G|\le  |R^*_G|\le 2n$. Then  $|B_G| + |R_G| \le 4n$. Finally, recall that $|B_G| \ge |A_{p-1}| \ge 3$ because   $A_{p-1} \subseteq B_G$.
We first consider the case when $ |R^*_G| \ge n$.   Since $|A_p|\ge n+1$, by Claim~\ref{n,n+1 Lemma}, $G[A_p]$ has no red edges.   
We   claim that $|B_G|\ge n$. Suppose $3 \le |B_G|\le n-1$.  Then $|A_p|\le (2n-1) \cdot 2^{k-2} + n-7$ by Claim~\ref{l-vertex}  applied to $B_G$ and $A_p$.  But then
  \[  
  |G|=|A_p|+|B_G|+|R^*_G|+|X\less X_r^*|\le [(2n-1) \cdot 2^{k-2} + n-7]+ (n-1)+2n+3(k-1)<n\cdot 2^k+1,
  \]
  for all  $k\ge3$,  a contradiction. Thus $|B_G|\ge n$, as claimed.   
  By  Claim~\ref{n,n+1 Lemma}, $G[A_p]$ has no blue edges $X_b^* = \es$, and so $|X''' |\le |X'' |\le k-1$.  Since $G[A_p\cup X']$ has neither red nor blue edges, it follows that $|A_p\cup X'|\le n\cdot 2^{k-2}$ by minimality of $k$.   But then
\[
|G|=|A_p\cup X'|+|X''\cup X'''|+(|B_G|+|R_G|)\le n\cdot 2^{k-2}+2(k-1)+4n<n\cdot 2^k+1,
\]
  for all $k\ge3$, a contradiction. 

It remains to consider the case when $|R^*_G|\le n-1$. If $|B^*_G| \ge  n-1$, by Claim \ref{l-vertex}(i,ii) applied to  $B^*_G$ and $A_p$,   we have 
\[ |A_p| + |B^*_G|+|X\less (X_r^*\cup X_b^*)|\le \begin{cases}  
(n \cdot 2^{k-1} +2) + (n-1)+3(k-2), &  \text{if } |B^*_G| = n-1\\
n \cdot 2^{k-1} + (2n+3)+3(k-2), &  \text{if } |B^*_G| \ge n.
 \end{cases}
 \]  
Thus in either case,  $ |A_p| + |B^*_G|+|X\less (X_r^*\cup X_b^*)|\le n \cdot 2^{k-1} + 2n+3k-3$.   But then 
\[
|G|=(|A_p|+|B^*_G|+|X\less (X_r^*\cup X_b^*)|)+|R^*_G|\le (n\cdot 2^{k-1}+ 2n+3k-3)+(n-1) <n\cdot 2^k+1, 
\]
for all $k \ge 3$, a contradiction.  Thus $ 3\le  |B^*_G| \le n-2$.  By Claim \ref{BR*_G}, either $|B_G^*| \ge4$ or $|R_G^*| \ge4$.   By applying Claim \ref{l-vertex}  to  $B_G^*$ when $|B_G^*| \ge4$ (or  $R_G^*$  when $|R_G^*| \ge4$) and $A_p$,  we have  $|A_p| \le (2n-3)\cdot 2^{k-1} + n-7$.  Then $|R^*_G|\ge n-2$, else  
\[
|G|=|A_p|+|B^*_G|+|R^*_G|+|X\less (X_r^*\cup X_b^*)|\le [(2n-3)\cdot 2^{k-1} + (n-7)]+(n-2) + (n-3)+3(k-2) <n\cdot 2^k+1, 
\]
for all $k \ge 3$ and $n\in\{6,7\}$, a contradiction. Thus  $n-2\le |R^*_G|\le n-1$.  By Claim \ref{l-vertex}(ii, iii) applied to $R_G^*$ and $A_p$, $|A_p| \le (21-2n) \cdot 2^{k-1-q} + (5n - 31)$. But then   
\begin{align*}
|G|&=|A_p|+|B^*_G|+|R^*_G|+|X\less (X_r^*\cup X_b^*)|\\
	&\le [(21-2n) \cdot 2^{k-1} +  (5n - 31)]+(n-2)  +  (n-1)+3(k-2) \\
		&<n\cdot 2^k+1,
\end{align*}

for all $k \ge 3$ and $n\in\{6,7\}$, a contradiction.  \medskip

This completes the proof of Theorem~\ref{C13C15}.\qed

\end{document}